\documentclass[final,a4paper]{siamltex}
\usepackage{amsmath,amssymb,latexsym}
\usepackage[final]{graphicx}
\usepackage{overpic}
\newcommand{\res}{\mathop{\rm res}}

\newcommand{\mult}{\mathop{\rm mult}}

\newcommand{\field}[1]{\mathbb{#1}}
\newcommand{\R}{\field{R}}
\newcommand{\Z}{\field{Z}}
\newcommand{\N}{\field{N}}
\newcommand{\C}{\field{C}}

\renewcommand{\AA}{{\mathcal A}}

\newcommand{\D}{\field{D}}

\newcommand{\FF}{{\mathcal S}}

\newcommand{\T}{{\field{T}}}

\renewcommand{\Re}{\mathop{\rm Re}}
\renewcommand{\Im}{\mathop{\rm Im}}

\newcommand{\isdef}{\stackrel{\text{\tiny def}}{=}}

\def\XXint#1#2#3{{\setbox0=\hbox{$#1{#2#3}{\int}$}
\vcenter{\hbox{$#2#3$}}\kern-.5\wd0}}

\title{Szeg\H{o} polynomials: a view from the Riemann-Hilbert window%
\thanks{This paper is dedicated to Ed Saff on the occasion of
his 60th birthday. }}%

\author{A.\ Mart\'{\i}nez-Finkelshtein\thanks{University of Almer\'{\i}a
and Instituto Carlos I de F\'{\i}sica Te\'{o}rica y Computacional, Granada
University, Spain.}
 }%
\date{This paper is dedicated to Ed Saff on the occasion of
}
\begin{document}
\maketitle

\begin{abstract}
This is an expanded version of the talk given at the conference
``Constructive Functions Tech-04''. We survey some recent results on
canonical representation and asymptotic behavior of polynomials
orthogonal on the unit circle with respect to an analytic weight.
These results are obtained using the steepest descent method based
on the Riemann-Hilbert characterization of these polynomials.

\end{abstract}



\pagestyle{myheadings} \thispagestyle{plain}
\markboth{MART\'{I}NEZ-FINKELSHTEIN, ANDREI}{SZEG\H{O} POLYNOMIALS: A
VIEW FROM THE RIEMANN-HILBERT WINDOW}

\section{Introduction}

During the Fall Semester of 2003 I was visiting the Department of
Mathematics of the Vanderbilt University, where I had the
opportunity to continue my collaboration with Ed Saff. I was very
excited with the evolution of the Riemann-Hilbert approach to the
asymptotic analysis of orthogonal polynomials, and discussed
extensively with Ed the new perspectives. He was the one who posed
the question: can this method tell us anything new about such a
classical object as the orthogonal polynomials on the unit circle
(OPUC, known also as Szeg\H{o} polynomials), in particular, about
their zeros? The question was more on a skeptical side. I was aware
of some previous work of the founders of the method,
\cite{MR2000e:05006}, \cite{deift99}, but none of these papers was
focused on the description of the zeros of the OPUC's. So, we
started to work and realized that we were able to find curious facts
even in the simplest situations. Later on I visited Ken McLaughlin,
at that time in Chapel Hill. A two-day discussion at Strong Caf\'{e} (a
recommended place) was crucial, and Ken joined the team. This paper
is a short and informal report on some of the advances we have had
so far.

Let me introduce some notation and describe the setting. For $r>0$,
denote  $\D_r \isdef \{z\in \C:\, |z|<r\}$, and  $\T_r \isdef \{z\in
\C:\, |z|=r \}$. A positive measure $\mu$ on $\T_1$ has the
Lebesgue-Radon-Nikodym decomposition
\begin{equation}\label{Radon}
d\mu(z)=w(z)\, |dz|+d\mu_s
\end{equation}
where $\mu_s$ is the singular part of $\mu$ with respect to the
Lebesgue measure on $\T_1$. Throughout, we will consider measures
satisfying the Szeg\H{o} condition
$$
\int_{\T_1} \log w(z) \, |dz|>-\infty\,,
$$
allowing to define the \emph{Szeg\H{o} function} (see e.g.\
\cite[Ch.\ X, \S 10.2]{szego:1975}):
\begin{equation}\label{standardSzego}
D (w; z) \isdef \exp\left( \frac{1}{4\pi }\,\int_0^{2\pi} \log
w(e^{i \theta}) \, \frac{e^{i \theta}+z}{e^{i \theta}-z}\, d\theta
\right)\,.
\end{equation}
This function is piecewise analytic and non-vanishing, defined for
$|z|\neq 1$, and we will denote by $D_{\rm i}$ and $D_{\rm e}$ its
values for $|z|<1$ and $|z|>1$, respectively, given by formula
\eqref{standardSzego}. It is easy to verify that
\begin{equation}\label{symmetry}
\overline{D_{\rm i}\left(w;\frac{1}{\overline{z}}
\right)}=\frac{1}{D_{\rm e}(w;z)}\,, \quad |z|>1\,,
\end{equation}
and for the boundary values we have
\begin{equation}\label{boundary_value_Szego}
w(z)=\frac{D_{\rm i}(w;z)}{D_{\rm e}(w;z)}=\frac{1}{|D_{\rm
e}(w;z)|^2}\,, \quad z\in \T_1\,.
\end{equation}
The first equality in \eqref{boundary_value_Szego} can be regarded
as a Wiener-Hopf factorization of the weight $w$, which is a key
fact for the forthcoming analysis.

For a nontrivial positive measure $\mu$ on $\T_1$ there exists a
unique sequence of polynomials $\varphi_n (z)=\kappa_n
z^n+\text{lower degree terms}$, $\kappa_n
>0$,  such that
\begin{equation}\label{orthogonalityConditions}
\oint_{\T_1} \varphi_n(z)  \overline{\varphi_m(z)}\, d\mu (z)
=\delta_{mn}\,, \quad m, n =0, 1, \dots
\end{equation}
We denote by $\Phi_n(z) \isdef  \varphi_n(z)/\kappa _n$ the
corresponding monic orthogonal polynomials. They satisfy the
Szeg\H{o} recurrence
$$
\Phi_{n+1}(z)=z \Phi_n(z)-\overline{\alpha _n}\, \Phi_n^*(z)\,,
\quad \Phi_0(z)\equiv 1\,,
$$
where we use the standard notation $\Phi_n^*(z)\isdef z^n
\overline{\Phi_n(1/\overline{z})}$. The parameters $\alpha
_n=-\overline{\Phi_{n+1}(0)}$ are called \emph{Verblunsky
coefficients} (also \emph{reflection coefficients} or \emph{Schur
parameters}) and satisfy $\alpha _n \in \D_1$ for $n=0, 1, 2,
\dots$. Furthermore, $\mu  \leftrightarrow \{\alpha _n\}$ is a
bijection; the map $\mu \rightarrow \{\alpha _n\}$ is an
\emph{inverse problem}, and it is known to be difficult. In
\cite{Simon05b} there is a thorough discussion of several techniques
to tackle this problem. Recently, the Riemann-Hilbert approach, not
described in \cite{Simon05b}, proved to be very promising in this
context also (see \cite{MR2000e:05006}, \cite{deift99},
\cite{deift05}, \cite{McLaughlin/Miller:2004}). The main goal of
this paper is a further discussion of how this method can shed new
light on the study of the asymptotics of the Szeg\H{o} polynomials.
We are not going to provide detailed proofs that can be found
elsewhere (the references are included), the aim is to show the
method in action and to discuss some new results in two apparently
simple situations.

The structure of the paper is as follows. Section \ref{sec:analytic}
is devoted to the case when $D(w;\cdot)$ is non vanishing and has an
analytic extension across $\T_1$. The main role here is played by
the \emph{scattering function}\footnote{Function $1/\FF$ is denoted
in \cite{Simon05a} by $r$, and in \cite[Section 6.2]{Simon05b} by
$b$. It corresponds also to the scattering matrix in
\cite{Geronimo/Case:79}.  I prefer to follow the notation of
\cite{math.CA/0502300}.}
\begin{equation}\label{def_F}
    \FF(w; z) \isdef D_{\rm i}(w; z)D_{\rm e}(w; z)\,,
\end{equation}
meromorphic in an annulus, containing $\T_1$, and which, via its
iterated Cauchy transforms, allows to write some canonical series
representing $\varphi_n$'s. In this situation convergence is always
exponentially fast, and the Riemann-Hilbert analysis is particularly
simple and transparent. Only some of the multiple corollaries of the
canonical representation are discussed; a more thorough analysis is
contained in \cite{math.CA/0502300}. In Section
\ref{sec:nonanalytic} we look at the situation when the original
analytic and nonvanishing weight has been modified by a factor
having a finite number of zeros on the unit circle. Now the behavior
of the zeros of $\Phi_n$'s is qualitatively different: most of them
cluster at $\T_1$, and only a finite number stays within the disc
$\D_1$. The method of Section \ref{sec:analytic} must be modified
now in order to handle the zeros of the weight: a local analysis
plays the major role. The exposition here is much more sketchy; in
this sense, more than a detailed view this window gives us a glimpse
of the possible techniques and results. At any rate, the main goal
is to persuade the reader that the Riemann-Hilbert analysis is a
powerfull technique, that deserves to be in the toolbox of those
interested in orthogonal polynomials on the unit circle.

\section{Analytic and nonvanishing weight}\label{sec:analytic}

Any analysis of OPUC's can be started either from the orthogonality
measure $\mu$ or from the sequence of the Verblunsky coefficients
$\{ \alpha _n\}$. Let us assume that the sequence of the Verblunsky
coefficients has an exponential decay:
\begin{equation}\label{resultNevaiTotik}
\rho\isdef \varlimsup_{n\to \infty} |\alpha _n|^{1/n}\,.
\end{equation}
Nevai-Totik \cite{Nevai/Totik:89} proved that this situation is
characterized by the following conditions on $\mu$: in the
decomposition \eqref{Radon}, $\mu_s=0$, measure $\mu$ satisfies the
Szeg\H{o} condition, and
\begin{equation}\label{Def_R_bis}
\rho =  \inf\{0<r<1:\, D_{\rm e}(w;z) \text{ is holomorphic in }
|z|>r \}\,.
\end{equation}
Taking into account \eqref{symmetry} we see that the first identity
in \eqref{boundary_value_Szego} can be regarded as an analytic
extension of the weight $w$. With this definition of $w$ (which we
use in the sequel) we can say equivalently that
\begin{equation}\label{Def_R}
\rho=\inf\{0<r<1:\, 1/w(z) \text{ is holomorphic in } r<|z|<1/r
\}\,,
\end{equation}
and both circles $\T_\rho $ and $\T_{1/\rho} $ contain singularities
of $1/w$. In this situation the well-known Szeg\H{o} asymptotic
formula
\begin{equation*}\label{ExteriorSzego}
\lim_{n\to \infty}\varphi_n^* (z) =\frac{1}{D_{\rm i}(w; z)}\,,
\quad z\in \D_1\,,
\end{equation*}
can be continued analytically through the unit circle $\T_1$ and is
valid locally uniformly in $\D_{1/\rho }$. It shows that the number
of zeros of $\{\varphi _n\}$ on compact subsets of $\D_1\setminus
\overline{\D}_\rho $ remains uniformly bounded, and these zeros are
attracted by the zeros of $D_{\rm e}$ in $\rho <|z|<1$
(\emph{Nevai-Totik points} in the terminology of B.\ Simon
\cite{Simon04b}). Numerical experiments show that the vast majority
of zeros gather at the ``critical circle'' $\T_\rho $. This fact was
justified theoretically by Mhaskar and Saff \cite{Mhaskar90}, who
using potential theory arguments showed that for any subsequence $\{
n_k \} \subset \N$ satisfying
$$
\rho=\lim |\alpha _{n_k}|^{1/n_k}\,,
$$
the zeros of $\{\varphi_{n_k+1}\}$ distribute asymptotically
uniformly in the weak-star sense on $\T_\rho $.

The behavior of the zeros inside $\D_\rho $ can be intriguing.
Although approaching in mass the critical circle $\T_\rho $ as
predicted by Mhaskar and Saff, some of them still may remain inside
and follow interesting patterns. Even the convergence to the circle
$\T_\rho $ is different for different measures. Can we give a full
description of this behavior in terms of the weight of
orthogonality? The answer is positive, and the description will
involve a sequence of iterates of certain Hankel and Toeplitz
operators with symbols depending on the scattering function $\FF$
introduced in \eqref{def_F}. It is a consequence of a canonical
representation of the Szeg\H{o} polynomials, found by means of the
Riemann-Hilbert characterization.

\subsection{Steepest descent analysis and canonical representation for orthogonal polynomials}
\label{sec:RH_analysis}

The starting point of all the analysis is the fact that under
assumption \eqref{resultNevaiTotik} conditions
\eqref{orthogonalityConditions} can be rewritten in terms of a
non-hermitian orthogonality for $\varphi_n$ and $\varphi_n^*$:
\begin{align*}
 \oint_{\T_1} \varphi_n(z)  z^{n-k-1} \,
\frac{w(z)}{z^n} dz & =0,\quad \text{for } k=0, 1, \dots, n-1\,,
 \\
 \oint_{\T_1} \varphi^*_{n-1}(z)  z^{k} \,
\frac{w(z)}{z^n} \, dz & = \begin{cases} 0, & k=0, 1, \dots, n-2, \\
i/\kappa _{n-1}, & k=n-1\,.
\end{cases}
\end{align*}
Here and in what follows, all the circles $\T_\alpha $, $\alpha >0$,
are oriented counterclockwise; with this orientation we talk about
the ``$+$'' and the ``$-$'' side of $\T_\alpha $ referring to its
inner and outer boundary points, respectively. Analogously, $f_+$
and $ f_- $ are the corresponding boundary values on $\T_\alpha $
for any function $f$ for which these limits exist. By standard
arguments (see e.g.\ \cite{MR2000e:05006} or \cite{MR2000g:47048},
as well as the seminal paper \cite{Fokas92} where the
Riemann-Hilbert approach to orthogonal polynomials started),
$$
Y(z)=\begin{pmatrix} \Phi_n(z) & \displaystyle \dfrac{1}{2\pi
i}\,\oint_{\strut\T_1}
\dfrac{\Phi_n(t) w(t)\, dt}{t^n(t-z)}  \\
-2\pi \kappa _{n-1} \varphi^*_{n-1}(z) & - \, \displaystyle
\dfrac{\kappa _{n-1}}{ i}\,\oint_{\strut\T_1}
\dfrac{\varphi^*_{n-1}(t) w(t)\, dt}{t^n(t-z)}
\end{pmatrix}
$$
is a unique solution of the following Riemann-Hilbert problem: $Y$
is holomorphic in $\C\setminus {\T_1}$,
\begin{equation}\label{RHproblem}
Y_+(t)=Y_-(t)\, \begin{pmatrix} 1 & w(t)/t^n \\ 0 & 1
\end{pmatrix}\,, \quad z\in \T_1\,, \quad \text{and} \quad
\lim_{z\to \infty} Y(z)\,\begin{pmatrix} z^{-n} & 0 \\ 0 & z^n
\end{pmatrix}=I\,,
\end{equation}
where $I$ is the $2 \times 2$ identity matrix.

This is the starting position for the steepest descent analysis as
described in \cite{MR2000g:47048} (see also \cite{Kuijlaars03}),
which consists in performing a series of \emph{explicit} and
\emph{reversible} steps in order to arrive at an equivalent problem,
which is solvable, at least in an asymptotic sense. Since these
steps are almost standard, they will
be described very schematically. We will use the following notation: $\sigma_3 =\begin{pmatrix} 1 & 0 \\
0& -1 \end{pmatrix}$ is the Pauli matrix, and for any non-zero
$x$ and integer $m$, $ x^{\sigma_3}=\begin{pmatrix} x & 0 \\
0& 1/x \end{pmatrix} $ and $x^{m\sigma _3}=(x^m)^{\sigma _3}$.

{\sc Step 1:} Define
\begin{equation}\label{defH}
H(z)\isdef \begin{cases} z^{-n \sigma _3} , & \text{if } |z|>1, \\
I, & \text{if } |z|<1,
\end{cases}
\end{equation}
and put $ T(z) \isdef Y(z)\, H(z) $. Then $T$ is holomorphic in
$\C\setminus {\T_1}$; this transformation normalizes the behavior at
infinity: $\lim_{z\to \infty} T(z)=I$. The price we pay is the
oscillatory behavior of the new jump matrix on $\T_1$:
$$
T_+(t)=T_-(t)\, \begin{pmatrix} t^n & w(t)  \\ 0 & t^{-n}
\end{pmatrix}\,, \quad t\in \T_1\,.
$$
We get rid of these oscillations in the next transformation, taking
advantage of the analyticity of its entries in the annulus.

\begin{figure}[htb]
\centering \begin{overpic}[scale=0.65]{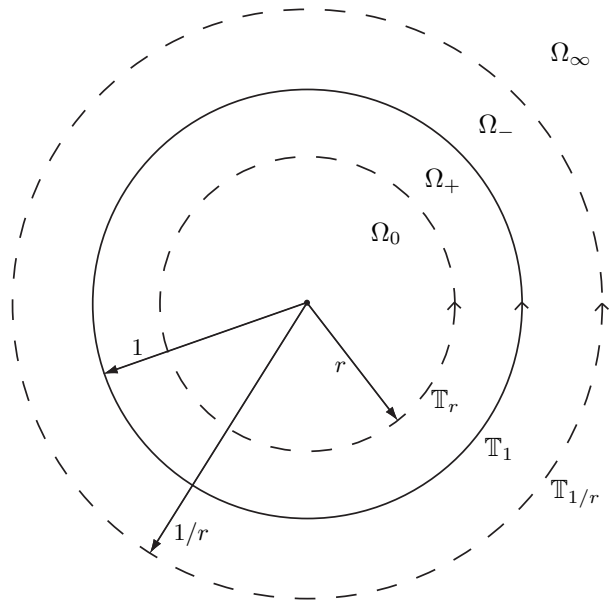}
      \put(79,24){$\T_1$}
         \put(70,32){$\T_r$}
          \put(90,17){$\T_{1/r}$}
          \put(60,60){$\Omega_0$}
          \put(69,69){$\Omega_+$}
          \put(78,78){$\Omega_-$}
          \put(90,90){$\Omega_\infty$}
          \put(54,38){\small $r$}
          \put(27,10){\small $1/r$}
          \put(20,41){\small $1$}
\end{overpic}
\caption{Opening lenses.}\label{fig:nice_case1}
\end{figure}

{\sc Step 2:} Choose an arbitrary $r$, $\rho <r<1$, that we fix for
what follows; it determines the regions (Figure
\ref{fig:nice_case1})
\begin{align*}
 \Omega_0 &= \{z:\, |z|<r \}, \qquad \Omega_\infty =\{z:\,
|z|>1/r \} ,\\
\Omega_+ &=\{z:\, r<|z|<1 \}, \qquad \Omega_- = \{z:\, 1<|z|<1/r \}
.
\end{align*}
Define $ U(z) \isdef T(z)K(z)$, where
\begin{equation}\label{KforT}
    K(z) \isdef \begin{cases} I, & \text{if } z \in \Omega_0 \cup \Omega_\infty, \\
\begin{pmatrix} 1 & 0  \\ z^n/  w(z)  & 1
\end{pmatrix}^{-1}\,, & \text{if }
z\in \Omega_+, \\ \begin{pmatrix}  1 & 0  \\ 1/(z^n w(z)) & 1
\end{pmatrix}\,, & \text{if }
z\in \Omega_-.
\end{cases}
\end{equation}
Then $U$ is holomorphic in $\C\setminus \left(  {\T_r} \cup {\T_1}
\cup {\T_{1/r}}\right)$. We have not modified its behavior at
infinity, but now
$$
U_+(t)=U_-(t)\, J_U(t), \quad t\in \left( {\T_r}  \cup {\T_1}\cup
{\T_{1/r}}\right),
$$
where
$$
J_U(t)=\begin{cases} \begin{pmatrix} 0 & w(t)  \\ -1/w(t) & 0
\end{pmatrix}, & \text{if } t\in {\T_1}, \\   \begin{pmatrix} 1 & 0  \\ t^n/  w(t)  & 1
\end{pmatrix}, & \text{if } t\in {\T_r}, \\
\begin{pmatrix} 1 & 0  \\ 1/(t^n w(t)) & 1
\end{pmatrix}, & \text{if } t\in {\T_{1/r}}.
\end{cases}
$$
The jump on $\T_r$ and $\T_{1/r}$ is exponentially close to the
identity, which is convenient to our purposes. We have to deal now
with the relevant jump on the unit circle.

{\sc Step 3:} The Szeg\H{o} functions $D_{\rm i}$ and $D_{\rm e}$
have been introduced in \eqref{standardSzego}. Define
 \begin{equation}\label{tau}
\tau  \isdef \frac{1}{D_{\rm i}(w; 0)}=D_{\rm e}(w;
\infty)=\exp\left( - \frac{1}{4\pi }\,\int_0^{2\pi} \log w(e^{i
\theta}) \, d\theta \right)
>0\,.
\end{equation}
Hence, if we introduce the geometric mean
\begin{equation}\label{def_G}
\mathcal G[w]\isdef   \exp\left(\frac{1}{2\pi} \int_{0}^{2\pi}
 \log\left( w\left(e^{i\theta} \right)\, d\theta \right)
 \right)\,,
\end{equation}
then $\tau =\left( \mathcal G[w]\right)^{-1/2}$.

It is straightforward to check that the piece-wise analytic
matrix-valued function
\begin{equation}\label{equ:defN}
    N(z)=N(w;z)\isdef \begin{cases}
    \begin{pmatrix} D_{\rm e}(w; z)/\tau  & 0  \\ 0 & \tau /D_{\rm e}(w; z)
\end{pmatrix}, & \text{if } |z|>1\,, \\
  \begin{pmatrix} 0 & D_{\rm i}(w; z)/\tau   \\ -\tau /D_{\rm i}(w; z) & 0
\end{pmatrix}, & \text{if } |z|<1\,,
    \end{cases}
\end{equation}
is invertible, and has the same jumps on $\T_1$ as $U(z)$. This
motivates to make a new transformation, defining $S(z) \isdef U(z)
N^{-1}(z)$. Matrix $S$ is holomorphic in $\C\setminus ({\T_r} \cup
{\T_1}\cup {\T_{1/r}} )$,
\begin{equation*}\label{boundsForJ_S}
\lim_{z\to \infty} S(z)=I,
\end{equation*}
and
\begin{equation}\label{jumpsForS}
S_+(t)=S_-(t)\, J_S(t), \quad t\in  {\T_r} \cup {\T_1}\cup
{\T_{1/r}},
\end{equation}
where
\begin{equation}\label{def_J_S}
J_S(t)=N_- J_U N_+^{-1}= \begin{cases} I, & \text{if } t\in {\T_1}, \\
\begin{pmatrix} 1 & -t^n \FF(w; t)/  \tau ^2  \\ 0  & 1
\end{pmatrix}, & \text{if } t\in {\T_r}, \\
\begin{pmatrix} 1 & 0  \\ \tau ^2/(t^n \FF(w; t)) & 1
\end{pmatrix}, & \text{if } t\in {\T_{1/r}}.
\end{cases}
\end{equation}
Our main character, $\FF$, has entered the picture!

Summarizing, we have
\begin{equation}\label{equ:exprForY}
    Y(z)=S(z)N(z)K^{-1}(z)H^{-1}(z)\,.
\end{equation}
Here $H$, $K$ and $N$ are explicitly defined in \eqref{defH},
\eqref{KforT} and \eqref{equ:defN}, respectively. About $S$ we know
only that it is piece-wise analytic and satisfies the jump condition
\eqref{jumpsForS}--\eqref{def_J_S}. A feature of this situation is
that we can write a formula for $S$ in terms of a series of iterates
of some Cauchy operators acting on the space of holomorphic
functions in $\C \setminus ({\T_r} \cup {\T_{1/r}} )$ with
continuous boundary values. Indeed, let us denote
$$
S=\begin{pmatrix}  S_{11} & S_{12} \\ S_{21} & S_{22}
\end{pmatrix}\,,
$$
and look at the equations given by the jumps on $\T_r$. According to
\eqref{def_J_S}, the first column is analytic across this circle,
while the second column has an additive jump equal to the first
column times $-t^n \FF(w; t)/  \tau ^2$. So, if we define the
operator
$$
\mathcal M_n^{\rm i}(f)(z)  \isdef -\frac{1}{2\pi i\, \tau^2}\,
\oint_{{\T_r}} f_-(t)\,
 \frac{\FF(w; t) \, t^n }{ t-z} \, dt\,,
$$
where $f_-$ denotes the exterior boundary values of the
function $f$ on $\T_r$, then taking into account the behavior at
infinity and Sokhotsky-Plemelj's theorem we get that
\begin{align}
\label{solvingForS1} S_{12}=\mathcal M_n^{\rm i}(S_{11})\,, \quad
 S_{22}=1+\mathcal M_n^{\rm
i}(S_{21})\,.
\end{align}
Analogously, with
$$
\mathcal M_n^{\rm e}(f)(z)   \isdef  \frac{\tau^2 }{2\pi i}\,
\oint_{{\T_{1/r}}} f_-(t)\, \frac{ 1 }{\FF(w; t)\, t^n\, (t-z)} \,
dt\,,
$$
where $f_-$ denotes the exterior boundary
values of the function $f$ on $\T_{1/r}$, we have again
\begin{align}
\label{solvingForS2} S_{11}=1+\mathcal M_n^{\rm e}(S_{12})\,, \quad
S_{21}=\mathcal M_n^{\rm e}(S_{22})\,.
\end{align}
By \eqref{solvingForS1}, \eqref{solvingForS2}, functions $S_{ij}$
satisfy the following integral equations:
\begin{align*}
(I-\mathcal M_n^{\rm e}\circ \mathcal M_n^{\rm i})S_{11} & =1\,,
\qquad (I- \mathcal M_n^{\rm i}\circ \mathcal M_n^{\rm e}) S_{12}
=\mathcal M_n^{\rm i}(1) \,, \\ (I- \mathcal M_n^{\rm e}\circ
\mathcal M_n^{\rm i}) S_{21} & =\mathcal M_n^{\rm e}(1)\,, \qquad
(I-\mathcal M_n^{\rm i}\circ \mathcal M_n^{\rm e}) S_{22} =1 \,,
\end{align*}
where $I$ is the identity operator. Straightforward bounds show that
there exists a constant $C>0$ depending on $r$ only, such that
\begin{equation}\label{bounds_operators}
\begin{split}
\left| \mathcal M_n^{\rm i}(f)(z) \right| & \leq C \, r^n\, \frac{
\|f_-\|_{\T_r}}{\left||z|-r \right|}\,, \qquad
z\notin \T_r\,, \\
 \left| \mathcal M_n^{\rm e}(f)(z) \right| & \leq
C \, r^n\, \frac{ \|f_-\|_{\T_{1/r}}}{\left||z|-1/r \right|}\,,
\qquad z\notin \T_{1/r}\,,
\end{split}
\end{equation}
where $\|\cdot \|_\gamma$ is the $\sup$-norm on $\gamma$; in the
sequel we use $C$ to denote some irrelevant constants, different in
each appearance, whose dependence or independence on the parameters
will be stated explicitly. Thus, we can invert these operators using
convergent Neumann series,
\begin{align*}
S_{11} & =\left( \sum_{k=0}^\infty \left(\mathcal M_n^{\rm e}\circ
\mathcal M_n^{\rm i}\right)^k \right)(1)\,, \qquad  S_{12} =\left(
\sum_{k=0}^\infty \left(\mathcal M_n^{\rm i}\circ \mathcal M_n^{\rm
e}\right)^k  \circ \mathcal M_n^{\rm i}\right)(1)  \,, \\ S_{21} &
=\left( \sum_{k=0}^\infty \left(\mathcal M_n^{\rm e}\circ \mathcal
M_n^{\rm i}\right)^k  \circ \mathcal M_n^{\rm e}\right)(1)\,, \qquad
S_{22} =\left( \sum_{k=0}^\infty \left(\mathcal M_n^{\rm i}\circ
\mathcal M_n^{\rm e}\right)^k \right)  (1) \,.
\end{align*}
In order to make this somewhat more explicit, let us define
recursively two sequence of functions:
\begin{equation*}\label{def_nf_1}
\begin{split}
 f_n^{(0)} & \isdef  1 \,, \quad  f_n^{(1)} \isdef  \mathcal M_n^{\rm
i}(1) \,, \quad \text{and} \quad f_n^{(2k)} \isdef  \mathcal
M_n^{\rm e} (f_n^{(2k-1)})\,, \quad f_n^{(2k+1)} \isdef  \mathcal
M_n^{\rm i} (f_n^{(2k)}) \,, \quad
 k \in
\N\,, \\
g_n^{(0)} & \isdef  1 \,, \quad  g_n^{(1)} \isdef  \mathcal M_n^{\rm
e}(1) \,, \quad \text{and} \quad g_n^{(2k)} \isdef  \mathcal
M_n^{\rm i} (g_n^{(2k-1)})\,, \quad g_n^{(2k+1)} \isdef  \mathcal
M_n^{\rm e} (g_n^{(2k)}) \,, \quad
 k \in \N\,.
 \end{split}
\end{equation*}
Then
\begin{equation*}\label{defSrepresentation}
\begin{split}
 S_{11} (n; z) &   =   \sum_{k=0}^\infty f_n^{(2k)}(z)\,, \quad S_{
21}(n; z)  = \sum_{k=0}^\infty g_n^{(2k+1)}(z)\,, \quad \text{ for }
|z|\neq 1/r\,, \\
  S_{12}(n; z) & = \sum_{k=0}^\infty f_n^{(2k+1)}(z)\,,
 \quad S_{22} (n; z)    =
\sum_{k=0}^\infty g_n^{(2k)}(z)\,, \quad \text{ for } |z|\neq r \,.
\end{split}
\end{equation*}
All these series are uniformly and absolutely convergent in their
domains. Moreover, straightforward bounds on the Cauchy transform
show that there exists a constant $C>0$ depending on $r$ only, such
that for $i=1, 2 $,
\begin{equation*}\label{uniformBoundsIandE}
|S_{i1} (z)|\leq \frac{C}{||z|-1/r|}\,, \text{ for } |z|\neq 1/r\,,
\quad \text{and} \quad |S_{i2} (z)|\leq \frac{C}{||z|-r|}\,, \text{
for } |z|\neq r\,.
\end{equation*}
We emphasize that in each of the regions above functions $S_{ij}$
have their own meaning, and are not obtained in general by analytic
continuation from one domain to another.

Once we have computed $S$, we may replace its expression in
\eqref{equ:exprForY} in order to find $Y$. It should be performed
independently in each region. For instance, in the domain
$\Omega_\infty$, $ N(z)= \left(  D_{\rm e}(w; z)/\tau
\right)^{\sigma _3}$, $K(z)=I$, $ H(z)=z^{-n\sigma _3}$, so that
$$
  Y(z)=S(z)(z^n D_{\rm e}(w; z)/\tau)^{\sigma _3}= \begin{pmatrix}  S_{11}(n; z)z^n D_{\rm e}(w;
  z)/\tau & S_{12}(n; z)\tau
  /(z^n D_{\rm e}(w; z)) \\ S_{21}(n; z)z^n D_{\rm e}(w; z)/\tau & S_{22}(n; z)\tau
  /(z^n D_{\rm e}(w; z))
\end{pmatrix}\,.
$$
Analogous computations are easily completed in the rest of the
regions.

All the information about the parameters of the orthogonal
polynomials is codified in the first column of $Y$: its $(1,1)$
entry gives us $\Phi_n$, that evaluated at $z=0$ yields the
Verblunsky coefficients, while  the $(2,1)$ entry at $z=0$ is
related to the leading coefficient $\kappa$: $Y_{21}(0)=-2\pi \kappa
_{n-1}^2$. Hence, we have obtained the following theorem, that has
been proved in \cite{math.CA/0502300}:
\begin{theorem}
\label{thm:nice case} Let $w$ be a strictly positive analytic weight
on the unit circle $\T_1$, the constant $\rho$ as defined in
\eqref{Def_R_bis}--\eqref{Def_R}, and constant $r$ with $\rho <r<1$
fixed. Then with the notations introduced above, for every $n\in \N$
sufficiently large the following formulas hold:
\begin{align}\label{representation}
i)& \qquad \Phi_n(z)=\begin{cases} \tau ^{-1} z^n D_{\rm e}(w; z)
\,S_{11}
(n; z)\,, & \text{if } |z|> 1/r\,; \\
\tau ^{-1} z^n D_{\rm e}(w; z) \,S_{11} (n; z)
-\dfrac{\tau\,S_{12}(n; z)}{ D_{\rm
 i}(w; z)}\,, & \text{if } r < |z|< 1/r\,; \\
 -\dfrac{\tau\,S_{12}(n; z)}{
D_{\rm
 i}(w; z)}\,, & \text{if }  |z|< r\,.
\end{cases}
\\
ii)& \qquad \overline{\alpha _n}=\tau^2 S_{12}(n+1; 0)\,.
\label{alphas}
\\
iii)& \qquad \kappa_{n}^2=\frac{\tau^2}{2 \pi}\,  S_{22}(n+1; 0) \,.
\label{kappas}
\end{align}
\end{theorem}

{\sc Remark:} It is easily seen that the method we have just
described is valid also if we replace the condition of positivity of
$w$ on $\T_1$ by the requirement that its winding number on $\T_1$
is zero. In such a case we can assure that $\deg \Phi_n=n$ only for
sufficiently large $n$'s, but the rest of the argument remains the
same.

\medskip

Before we analyze some implications of these formulas let us look
more carefully at the scattering function $\FF$ and at the iterates
of its Cauchy transforms used in the definition of $S$. Following
notation of \cite[Section 6.2]{Simon05b}, let
$$
\log w(z)=\sum_{k\in \Z} \hat L_k z^k
$$
be the Laurent expansion (equivalently, the Fourier series) for
$\log w$. Then straightforward computation shows that
$$
\FF(w; z)=\exp\left( \sum_{k=1}^{+\infty} \left(\hat L_k z^k - \hat
L_{-k} z^{-k} \right)\right)=\exp\left( \sum_{k=1}^{+\infty}
\left(\hat L_k z^k - \overline{\hat L_{k}} z^{-k} \right)\right)\,.
$$
Let
\begin{equation}\label{LaurentExpansion}
\FF(w; z)=\sum_{k=-\infty}^{+\infty } \left( \FF \right)_k z^k \quad
\text{and} \quad \frac{1}{\FF(w; z)}=\sum_{k=-\infty}^{+\infty }
\left( \frac{1}{\FF} \right)_k z^k
\end{equation}
be the Laurent expansions of $\FF$ and $1/\FF$ in the annulus $\rho
<|z| <1/\rho $, respectively. Taking into account that $|\FF(w;
z)|=1$ on $\T_1$,
\begin{equation}\label{LaurentFminus1}
\left( \frac{1}{\FF}  \right)_k= \overline{\left( \FF
\right)_{-k}}\,, \quad \ \text{and} \quad \sum_{k=-\infty}^{+\infty
} \left( \FF \right)_{k+m} \overline{\left( \FF
\right)_k}=\begin{cases}  1, & \text{if } m=0, \\ 0, & \text{if } m
\in \Z \setminus \{ 0\}\,.
\end{cases}
\end{equation}
Denote by $H^2_+$ the Hardy class, and $H_-^2\isdef L^2 \ominus
H^2_+$. Let us denote by $\mathcal P_+$ and $\mathcal P_-$ the Riesz
projections onto $H^2_+$ and $H^2_-$, respectively, and $\sigma_n
(z) \isdef z^n \FF(w; z) $. Then
$$
\mathcal P_+ (\sigma_n)(z)=\sum_{k=-n}^{+\infty } \left( \FF
\right)_k z^{k+n}\,, \qquad \mathcal P_-(\sigma_n)(z)=\sum_{k<-n}
\left( \FF \right)_k z^{k+n}\,.
$$
In particular,
\begin{equation*}\label{f1_Plus}
f_n^{(1)}(z)= \mathcal M_n^{\rm i}(1) =\begin{cases}
 \displaystyle -\frac{1}{\tau ^2}\,  \mathcal P_+  (\sigma_n )(z)
 =-\frac{1}{\tau ^2}\,  \sum_{k=-n}^{+\infty } \left( \FF
\right)_k z^{k+n} \in H_+^2\,, & \text{if } |z|< r \,, \\
\displaystyle  \frac{1}{\tau ^2}\,  \mathcal P_- (\sigma_n )(z) =
\frac{1}{\tau ^2}\, \sum_{k<-n} \left( \FF \right)_k z^{k+n} \in
H^2_-\,, & \text{if } |z|> r\,.
 \end{cases}
\end{equation*}
The series in the right hand side converge locally uniformly.
Observe also that
\begin{equation}\label{connectionfandc}
f_n^{(1)}(0)=- \frac{1}{\tau^{2}}\, \left( \FF \right)_{-n}=-
\frac{1}{\tau^{2}}\, \overline{\left( \frac{1}{\FF} \right)_{n}}\,.
\end{equation}

If we introduce the following composition of Hankel and Toeplitz
operators, having $\sigma _n^{\pm 1}$ as a symbol,
$$
\mathcal H_n^\pm:\, H^2_- \mapsto H^2_\pm\,, \quad \text{given by }
\quad \mathcal H_n^+(f)=-  \mathcal P_+(\sigma _n \mathcal P_-
(\sigma _n^{-1} f ) ) \,, \quad
 \mathcal H_n^-(f)=  \mathcal P_-(\sigma _n \mathcal P_-
(\sigma _n^{-1} f ) ) \,,
$$
then
$$
f_n^{(2k+1)}=\begin{cases} \mathcal H_n^+(f_{n,-}^{(2k-1)})\,, &
\text{if } |z|< r \,, \\ \mathcal H_n^-(f_{n,-}^{(2k-1)})\,, &
\text{if } |z|> r \,,
\end{cases}
$$
where $f_{n,-}^{(2k-1)}$ represents the values of $f_{n}^{(2k-1)}$
in $\C\setminus \overline{\D_r}$.

For $g_n^{(k)}$ we can obtain analogous formulas:
$$
g_n^{(2)}(z)= \mathcal M_n^{\rm i} (g_n^{(1)})= \begin{cases}
 \displaystyle -\tau ^{-2}\, \mathcal P_+(\sigma _n  g_n^{(1)} )(z)=
   -\mathcal P_+(\sigma _n \mathcal P_+(\sigma _n^{-1} ) )(z) \,, & \text{if } |z|< r \,, \\
\displaystyle   \tau ^{-2}\, \mathcal P_-(\sigma _n  g_n^{(1)} )(z)=
\mathcal P_-(\sigma _n \mathcal P_+(\sigma _n^{-1} ) )(z) \,, &
\text{if } |z|> r\,,
 \end{cases}
$$
and
$$
g_n^{(2k+2)}=\begin{cases} \mathcal H_n^+(g_{n,-}^{(2k)})\,, &
\text{if } |z|< r \,, \\ \mathcal H_n^-(g_{n,-}^{(2k)})\,, &
\text{if } |z|> r \,.
\end{cases}
$$
For instance, taking into account \eqref{LaurentFminus1},  for
$|z|<r$,
\begin{equation*}
\begin{split}
 g_n^{(2)}(z) & = -\mathcal P_+(\sigma _n \mathcal P_+(\sigma _n^{-1} )
)(z)=-\mathcal P_+\left(\left(\sum_{j=-\infty}^{+\infty } \left( \FF
\right)_j z^{j+n}\right) \left(\sum_{k= n }^{+\infty } \left(
\frac{1}{\FF} \right)_{k} \,
z^{k-n}\right) \right) \\
&= -  \sum_{k\leq j, \;k\leq -n } \left( \FF \right)_j
\overline{\left( \FF \right)_k}\, z^{j-k}=- \sum_{m=0}^{+\infty}
\left( \sum_{k\leq -n } \left( \FF \right)_{k+m} \overline{\left(
\FF \right)_k}\right)\, z^{m} \,,
\end{split}
\end{equation*}
and by identities in \eqref{LaurentFminus1} we can rewrite last
formula as
\begin{equation} \label{g1}
1+ g_n^{(2)}(z) = \sum_{m=0}^{+\infty} \left( \sum_{k> -n } \left(
\FF \right)_{k+m} \overline{\left( \FF \right)_k}\right)  z^{m} \,,
\quad |z|<r\,.
\end{equation}

{\sc Remarks:} There are some further relations and equivalent
expressions that an interested reader can easily derive. For
instance, if we introduce the operator $\mathcal T_n$ on $L^2$ with
kernel
$$
T_n(i,j)\isdef \sum_{s=0}^{+\infty} \left( \FF \right)_{-(s+n+1+j)}
\left( \frac{1}{\FF}  \right)_{s+n+1+i}=\sum_{k<-n-j} \left( \FF
\right)_{k} \left( \frac{1}{\FF}  \right)_{i-j-k}\,,
$$
then
$$
\mathcal M_n^{\rm e}\circ \mathcal M_n^{\rm i}(f)=\begin{cases}
\mathcal P_+ \mathcal T_n (f) \,, & \text{if } |z|<1/r\,,  \\
\mathcal -P_- \mathcal T_n (f) \,, & \text{if } |z|>1/r\,,
\end{cases}
$$
and we can obtain expressions for $f_n^{(2k)}$ and $ g_n^{(2k+1)}$.

\subsection{Asymptotic behavior of OPUC} \label{sec:asymptotics}

Representation \eqref{representation} is asymptotic in nature. Using
\eqref{bounds_operators}, it is immediate to show that for all
sufficiently large $n$ and for $N=0, 1, 2,\dots $,
\begin{equation*}\label{errorTruncation1}
    \left| S_{11} (n; z)  -   \sum_{k=0}^N f_n^{(2k)}(z)
     \right| \leq \frac{C }{\left| |z|-1/r\right|}\,  r^{(2N+2)
    n}\,, \quad |z|\neq 1/r\,,
\end{equation*}
where the constant $C $ depends only on $r$ and $N$, but neither on
$n$ nor on $z$. Analogously,
\begin{equation}\label{errorTruncation2}
    \left| S_{12} (n;z)  - \sum_{k=0}^N f_n^{(2k+1)}(z)
    \right| \leq \frac{C }{\left| |z|-r\right|}\,  r^{(2N+3)
    n}\,, \quad |z|\neq r\,,
\end{equation}
where $C $ has a similar meaning as above. These bounds show that
\eqref{representation} allows us to obtain approximations of
$\{\Phi_n\}$ of an arbitrarily high order. We will concentrate only
on the most interesting domain including the critical circle
$\T_\rho $ and its interior (for a full analysis, check
\cite{math.CA/0502300}).

Let us discuss the consequences of truncating $S_{11}$ and $S_{12}$
in \eqref{representation} at their first terms,
$$
S_{11}(n; z)= 1+\mathcal O(r^{2n})\,, \quad S_{12}(n;
z)=-\frac{1}{2\pi i \tau ^2}\, \oint_{{\T_r}}
 \frac{\FF(w; t) \, t^n }{ t-z} \, dt+\mathcal O(r^{3n})\,,
$$
imposing some additional conditions on the analytic continuation of
our weight $w$ (or function $\FF$). We assume first that the
critical circle $\T_\rho $ contains only isolated singularities in a
finite number.
\begin{theorem}[\cite{math.CA/0502300}] \label{prop:residues}
Assume that there exists $0\leq \rho '<\rho $ such that $D_{\rm e}$
can be continued to the exterior of the circle $\T_{\rho'}$, as an
analytic function whose only singularities are on the circle
$\T_\rho$, and these are all isolated. Denote by $a_1, \dots, a_u$
the singularities (whose number is finite) of $D_{\rm e}$ on
$\T_\rho $,
\begin{equation*}\label{equal_poles}
|a_1|= \dots= |a_u|=\rho \,.
\end{equation*}
Then for $\rho <r'<r$ there exist constants $0\leq \delta=\delta
(r')<1$ and $C=C(r')<+\infty$ such that for $\rho '<|z|\leq r'$ and
$n\in \N$,
\begin{equation}\label{case_finite_number_sing2}
\left|\Phi_n(z)-  z^n\,  \frac{D_{\rm e}(w; z)}{D_{\rm e}(w; \infty)
} -\frac{ D_{\rm i}(w; 0)}{  D_{\rm i}(w; z)} \sum_{k=1 }^u
\res_{t=a_k } \left( \FF (w; t) \, \frac{t^n }{ t-z } \right)\right|
\leq C  \left( \rho ^n \delta ^n+ r^{3n}\right)\,.
\end{equation}
Furthermore, for every compact set $K\subset\D_\rho$ there exist
constants $0\leq \delta=\delta (K)<1$ and $C=C(K)<+\infty$ such that
for $z \in K$,
\begin{equation}\label{asymptotics_finite_sing}
\left| \Phi_n(z)-  \frac{D_{\rm i}(w; 0)}{  D_{\rm i}(w; z)}
\sum_{k=1 }^u \res_{t=a_k } \left( \FF (w; t)   \, \frac{t^n }{ t-z
} \right)\right| \leq C \left( \rho ^n \delta ^n+ r ^{3n} \right)\,.
\end{equation}
\end{theorem}

If $D_{\rm e}$ can be continued as an analytic function with a
finite number of isolated singularities to whole disc $\D_1$, then
we may take $\delta =0$ in the right hand sides in
\eqref{case_finite_number_sing2}--\eqref{asymptotics_finite_sing}.
Otherwise the right hand side in \eqref{asymptotics_finite_sing} may
be replaced by an estimate of the form $C \rho ^n \delta ^n$.

In order to isolate the zeros of $\Phi_{n}$, one must be able to
analyze the approximation to $\Phi_{n}$ afforded by
\eqref{case_finite_number_sing2} and
\eqref{asymptotics_finite_sing}. For example, zero-free regions may
be determined by \emph{(i)} establishing zero-free regions for the
approximation, and \emph{(ii)} bounding $\Phi_{n}$ away from zero
using the error estimates. Similarly, isolating the zeros can be
done by first isolating the zeros of the approximation, and then
using a Rouche' type argument for $\Phi_{n}$.

This theorem tells us that in general all the relevant information
for the asymptotics of $\Phi_n$'s in $\D_r$ comes from the
singularities of the exterior Szeg\H{o} function $D_{\rm e}$ on
$\T_\rho $ (that is, from the first singularities of $D_{\rm e}$ we
meet continuing it analytically inside the unit disc), and reduces
the asymptotic analysis of $\Phi_n$'s (at least, in the first
approximation) to the study of the behavior of the corresponding
residues. In the case when all the singular points that we met on
$\T_\rho $ are poles, this analysis is more or less straightforward.
\begin{definition}
Let $a\in \D_1$ be a pole of a function $f$ analytic in $|z|<1$. We
denote by $ \mult_{z=a}f(z)$ its multiplicity and say that $a$ is a
\emph{dominant} pole of $f$ if for any other singularity $b$ of $f$,
either $|a|>|b|$ or $|a|=|b|$, but then $b$ is also a pole and
$$
\mult_{z=a}f(z) \geq \mult_{z=b}f(z)\,.
$$
\end{definition}
In the sequel we use the following notation: for $a\in \C$ and
$\varepsilon >0$,
\begin{equation}\label{def_neighborhoods}
B_\varepsilon (a)\isdef \{z\in \C:\, |z-a|<\varepsilon  \}\,.
\end{equation}

\begin{theorem}[\cite{math.CA/0502300}] \label{prop:many dominant}
Assume that there exists \/ $0\leq \rho '<\rho $ such that $D_{\rm
e}$ can be continued to the exterior of the circle $\T_{\rho'}$, as
a meromorphic function whose only singularities are on the circle
$T_\rho$. Denote by $a_1, \dots, a_u$ the poles (whose number is
finite) of $D_{\rm e}$ on $\T_\rho $, and assume that the dominant
poles of $D_{\rm e}$ are $a_1, \dots, a_\ell$, $\ell \leq u$, and
their multiplicity is $m$.

Let $\varepsilon >0$. Then for $\rho '<|z|\leq r-\varepsilon$, $z
\notin \cup_{k=1}^u B_\varepsilon (a_k)$, and $n\in \N$,
\begin{equation}\label{case_finite_number_dominant_sing}
\Phi_n(z)=  \frac{D_{\rm e}(w; z) }{D_{\rm e}(w; \infty)}\, z^n
+\frac{ D_{\rm i}(w; 0)}{  D_{\rm i}(w; z)} \, \sum_{k=1 }^\ell
\binom{n}{m-1}\, a_k^{n-m+1}\, \frac{D_{\rm i} (w; a_k)\, \widehat
D_{\rm e} (w; a_k ) }{ a_k-z
 } \,   + h_n(z)   \,,
\end{equation}
where $ \widehat D_{\rm e} (w; a_k)=\lim_{z\to a_k} D_{\rm e} (w;
z)(z-a_k)^m $, $ k=1, \dots, \ell$. There exist a constant
$0<C<+\infty$ independent of $n$ and $\varepsilon $, and a constant
$0<\delta=\delta (\varepsilon )<1$, such that
\begin{equation*}\label{newBOundForH}
    |h_n(z)|\leq \begin{cases}
    C \left( \rho ^n \delta ^n+ r^{3n}\right), & \text{if } \, m=1\,, \\
    \dfrac{C}{ \varepsilon ^{m-1}}\, n^{m-2} \, \rho ^n\,, & \text{if }\, m\geq 2\,.
    \end{cases}
\end{equation*}

Furthermore, for every compact set $K\subset\D_\rho$ there exists a
constant $C=C(K) <\infty$ such that for $z \in K$, and $n\in \N$,
\begin{equation}\label{asymptotics_dominant_poles}
\left| \frac{\tau \, D_{\rm i}(w; z)}{a_1 ^{n-m+1}}\, \binom{n }{
m-1 }^{-1}\,\Phi_n(z)-
   \sum_{k=1 }^\ell \frac{D_{\rm i} (w; a_k)\,  \widehat D_{\rm e} (w; a_k) }{  a_k-z
 }   \, e^{2\pi i (n-m+1) \theta_k}\right|\leq \begin{cases}
C   \delta ^n, & \text{if } \, m=1\,, \\
    \dfrac{C}{   n} \,, & \text{if }\, m\geq 2\,,
 \end{cases}
\end{equation}
where
\begin{equation}\label{defTheta}
\theta_1=1, \quad \text{and} \quad
\theta_k=\frac{1}{2\pi}\,\left(\arg a_k - \arg a_1 \right),
  \qquad  k=2, \dots, \ell\, .
\end{equation}
In particular, on every compact set $K\subset \D_\rho $, for all
sufficiently large $n$ polynomials $\Phi_n$ can have at most $\ell
-1$ zeros, counting their multiplicities.
\end{theorem}

Observe that this result is applicable to weights of the form
$w(z)=|R(z) S(z)|^2$, $z\in \T_1$, where $R$ is a rational function
with at least one zero on $\T_\rho $ (or one pole on $\T_{1/\rho
}$), and $S$ is any function holomorphic and $\neq 0$ in any
annulus, containing $\{ \rho \leq |z| \leq 1/\rho \}$.

By means of \eqref{case_finite_number_dominant_sing} we may show
also that under assumptions of Theorem \ref{prop:many dominant} the
vast majority of zeros approaching the critical circle $\T_\rho $
does it in an organized way, exhibiting an equidistribution pattern:
if zeros $z^{(n)}_j$ of $\Phi_n$ can be numbered in such a way that,
roughly speaking,
\begin{equation*}\label{asymptModEquidistr}
|z^{(n)}_i|=\rho \left(1+\frac{1}{n}\, \log \binom{n}{m-1}+\mathcal
O\left(\frac{1}{n}\right) \right)\,,
\end{equation*}
and
\begin{equation*}\label{asymptArgEquidistr}
\arg\left(z^{(n)}_{i+j}\right)-\arg\left(z^{(n)}_i\right)=\frac{2\pi
j}{n}+\mathcal O \left( \frac{1}{n^2} \right)\,
\end{equation*}
(again, we refer the reader to \cite{math.CA/0502300} for details).
This is also an analogue of the interlacing property of the zeros of
orthogonal polynomials on the real line. Moreover,
\eqref{asymptotics_dominant_poles} allows us to describe the
accumulation set of zeros of $\{ \Phi_n\}$'s inside $\D_\rho $. For
instance, if all $\theta_k$'s in \eqref{defTheta} are rational, this
set is discrete and finite. Otherwise, as it follows from
Kronecker-Weyl theorem, it can be a diameter of $\D_\rho $ or even
fill a two dimensional domain.

\begin{figure}[htb]
\centering
\begin{tabular}{ll}
\hspace{-1.5cm}\mbox{\includegraphics[scale=0.65]{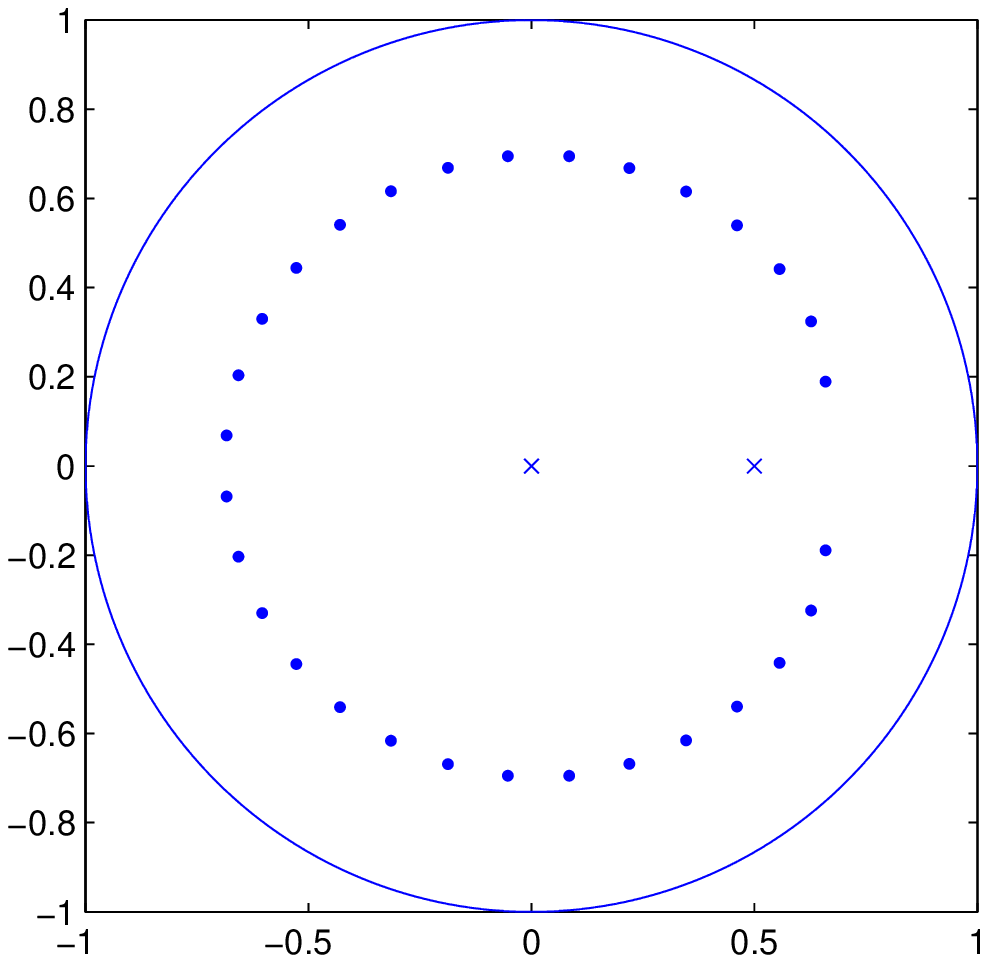}}
&
\hspace{-2cm}\mbox{\includegraphics[scale=0.65]{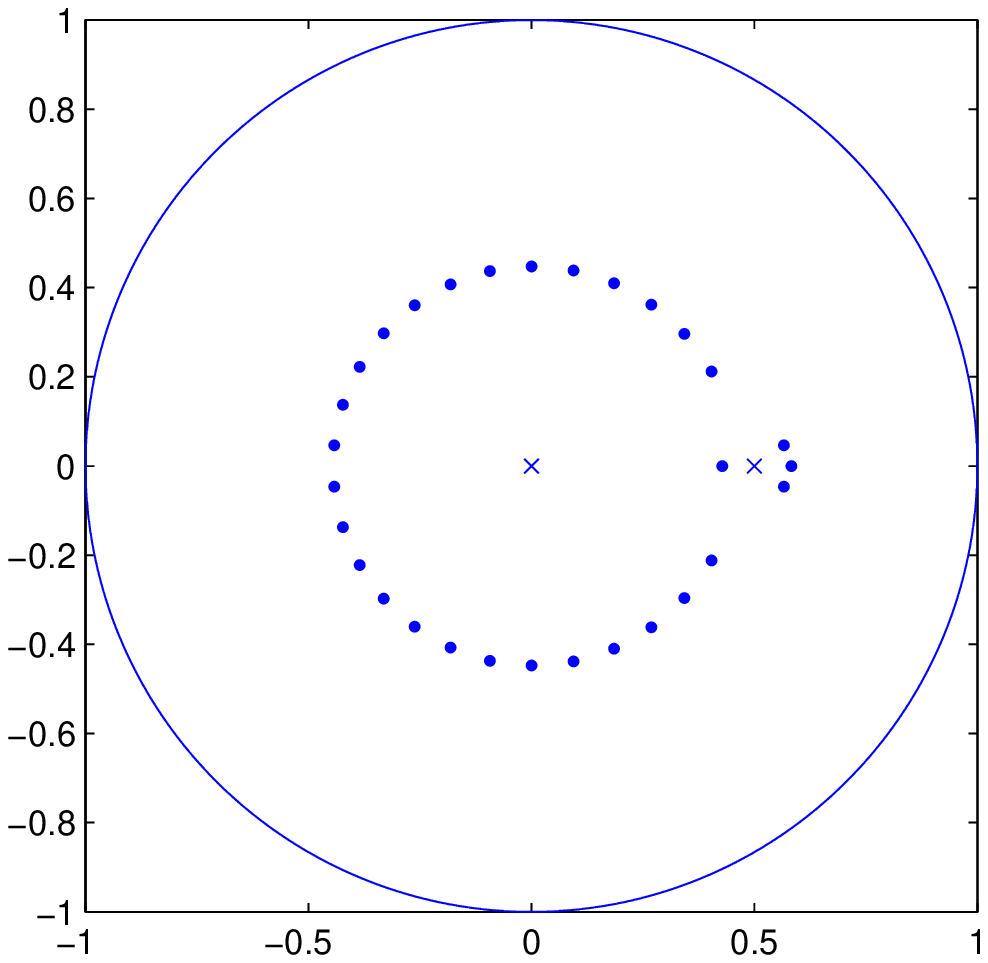}}
\end{tabular}
 \caption{Zeros of $\Phi_{30}$ for weights $w$ (left) and $1/w$ (right), with $w$
given in \eqref{weightExampleSing} and $\rho =1/2$.
}\label{fig:Essential30_dominantExample}
\end{figure}

The situation gets much more complicated if the first singularity
that we meet continuing $D_{\rm e}$ analytically inside is more
severe. Consider the simplest example of an essential singularity on
$\T_\rho $:
\begin{equation}\label{weightExampleSing}
w(t)= \left|\exp\left( \frac{1}{\rho -t}\right) \right|^2\,, \quad t
\in \T_1\,,
\end{equation}
with $0<\rho <1$. Observe that its inverse, $1/w$, satisfies also
the conditions of Theorem \ref{prop:residues}. However, the behavior
of the zeros of the OPUC for $w$ and $1/w$ is qualitatively
different, check Figure \ref{fig:Essential30_dominantExample}.

\begin{figure}[htb]
\centering
\begin{tabular}{ll}
\hspace{-1.5cm}\mbox{\includegraphics[scale=0.65]{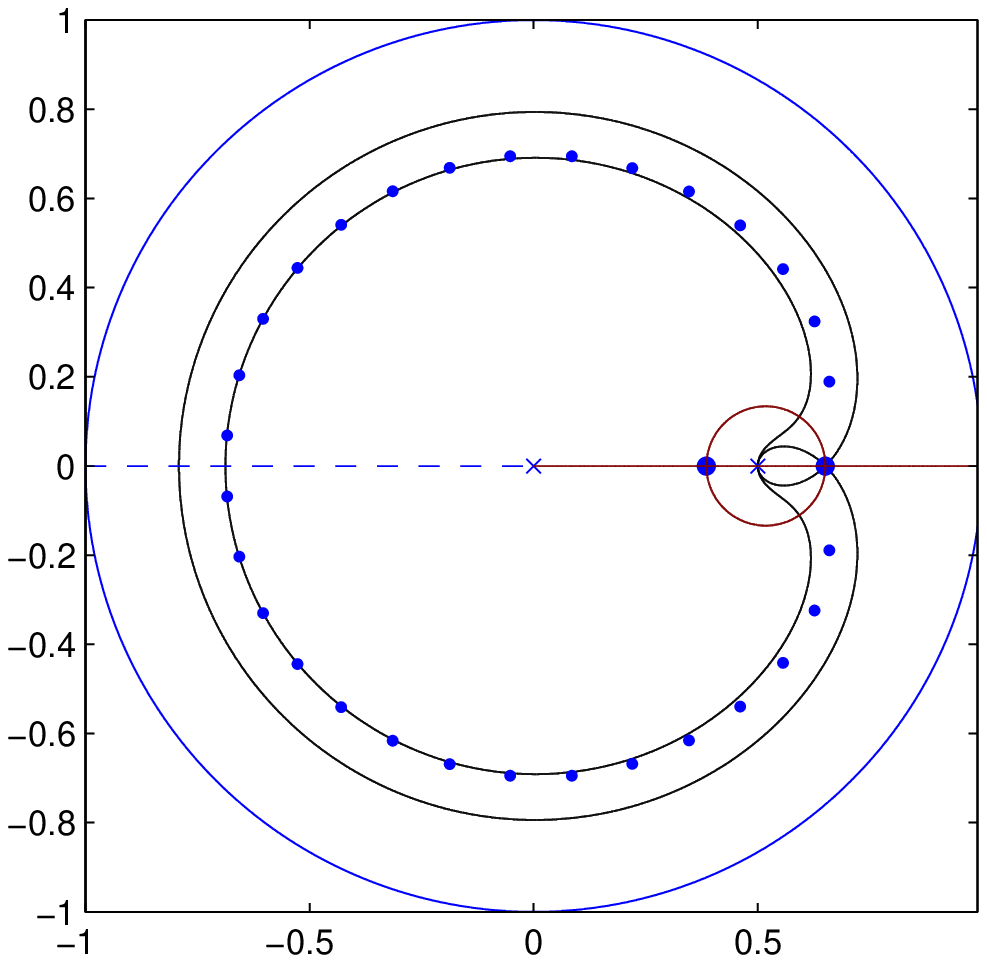}}
&
\hspace{-2cm}\mbox{\begin{overpic}[scale=0.65]{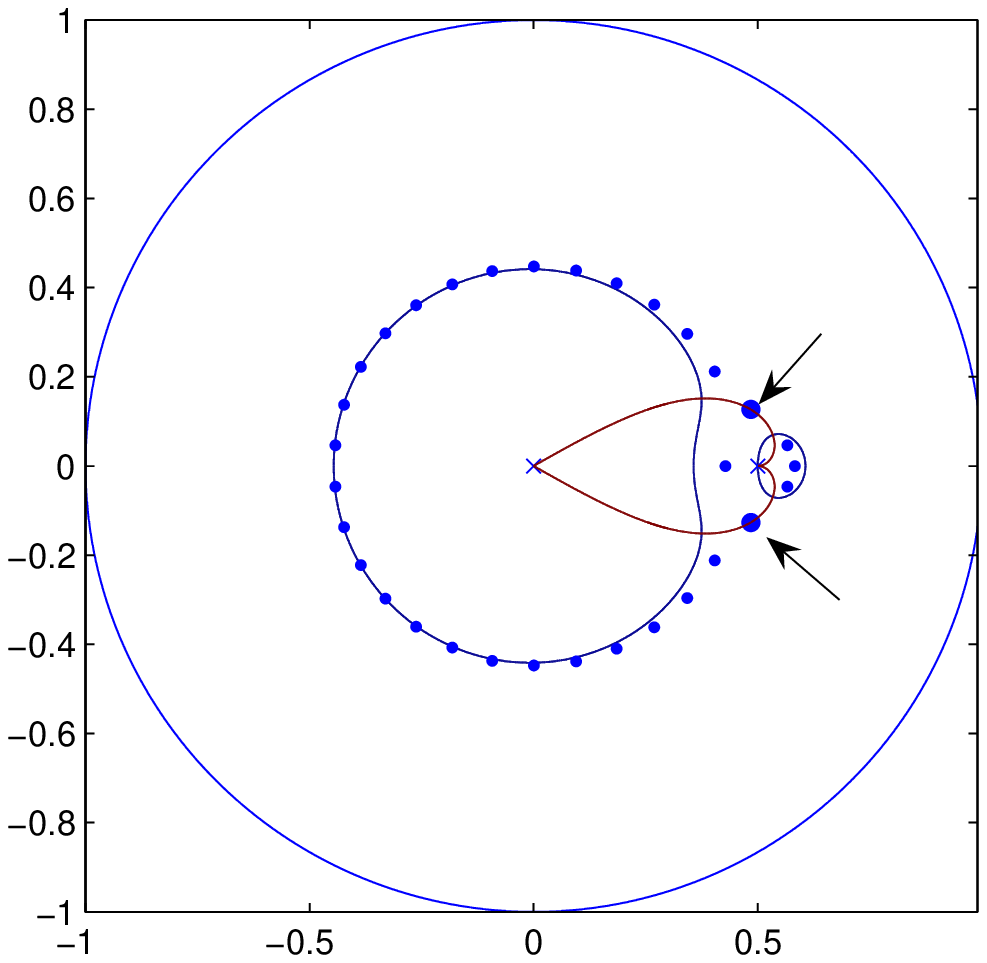}
\put(72,49){\small $t_+$}
    \put(74,27){\small $t_-$}
   \put(49,37){\small $0$}
\end{overpic}}
\end{tabular}
 \caption{Left: zeros of $\Phi_{n}$ for $w$ given in \eqref{weightExampleSing} with
$\rho=1/2$ and $n=30$, along with the level curves $\Re (\Psi _n(z)-
\Psi _n(t_+))=0$, $\Re ( \Psi _n(z)- \Psi _n(t_+))=\frac{1}{n}\,
\log\left(\frac{1}{2\sqrt{\pi} } \frac{\rho^{3/4}}{ n^{3/4}}
\right)$ and $\Im (\Psi _n(z)- \Psi _n(t_+))=0$. Right: the same,
but for the weight $1/w$.  Observe that the level curve
\eqref{levelCurve} has now two
components.}\label{fig:Essential30_dominant}
\end{figure}

In a few words, the explanation for this phenomenon is the
following. Observe that in the case of an essential singularity of
$\FF$ the asymptotic behavior of the Cauchy transform
$$
f_n^{(1)}(z)=-\frac{1}{2\pi i \tau ^2}\, \oint_{\T_r} \frac{t^n
\FF(w; t)}{t-z}\, dt
$$
is not as simple as when the only singular points are poles. In
fact, the leading term of the asymptotics will come now from a
dominant saddle point of
\begin{align*}
\Psi _n(t) &\isdef  \log t + \frac{1}{n}\, \log \FF(w; z)
\end{align*}
lying close to the singular point of $\FF$, $t=t_+$, which is the
solution of the equation
\begin{equation}\label{equationForT}
\frac{1}{t}-\frac{1}{n+1} \, \left(\frac{t}{\rho
t-1}+\frac{1}{t-\rho } \right)=0 \quad \text{satisfying} \quad t_+=
\rho  + \sqrt{\frac{\rho }{n+1}}+\mathcal
O\left(\frac{1}{n}\right)\,, \quad n \to
    \infty\,,
\end{equation}
where we take the positive square root. It is possible to show that
the zeros of the orthogonal polynomials (at least those not staying
to close to $z=\rho $), will approach the level curve
\begin{equation}\label{levelCurve}
\Re ( \Psi _n(z)- \Psi _n(t_+))=\frac{1}{n}\,
\log\left(\frac{1}{2\sqrt{\pi} } \frac{\rho^{3/4}}{ n^{3/4}}
\right)\,,
\end{equation}
and the error decreases with $n$. However, for the weight $1/w$ we
will have two dominant saddle points, and the different structure of
the level curve \eqref{levelCurve} for weights $w$ and $1/w$
explains the different result of the numerical experiments (see
Figure \eqref{fig:Essential30_dominant}).

\subsection{Verblunsky and leading coefficients}

Let us finish this Section with some comments about the behavior of
the notorious coefficients related to the OPUC.

Evaluating the polynomial $\Phi_{n+1}$ or any of its approximations
at the origin we obtain information about the Verblunsky
coefficients. For instance, a combination of \eqref{alphas},
\eqref{connectionfandc}, and  \eqref{errorTruncation2} yields the
following estimate of the Verblunsky coefficients $\alpha _n$:
\begin{proposition}[\cite{math.CA/0502300}]\label{cor:verblunsky2}
Let $w$ be a strictly positive analytic weight on the unit circle
$\T_1$. With the notation introduced above and for each $n\in \N$,
\begin{equation}\label{asympt_expr_alpha}
\alpha _n=-\left( \frac{1}{\FF}\right)_{n+1}+\mathcal O (r^{3n}) \,,
\end{equation}
where $(1/\FF)_{n+1}$ is the corresponding Laurent coefficient of
$1/\FF$ in \eqref{LaurentExpansion}.
\end{proposition}

This fact has the following reading: consider the generating
function of the Verblunsky coefficients,
$$
G(z)=\sum_{n=0}^\infty \alpha _n z^n\,.
$$
Then the Maclaurin series of $G$ and $\mathcal P_+(-z/\FF(z))$ match
up to the $\mathcal O (r^{3n})$ term. In consequence, we have
\begin{proposition}\label{prop:generating}
Function
$$
G(z)+\frac{z}{\FF(w; z)}\,,
$$
defined in a neighborhood of $\T_1$, can be continued as a
holomorphic function to the annulus $1<|z|<1/\rho^3$.
\end{proposition}
This fact has been established independently by Simon
\cite{Simon05a} and Deift and \"{O}stensson \cite{deift05}.

If the only singularities on $\T_\rho$ are dominant poles, we can
use formula \eqref{asymptotics_dominant_poles} in order to derive
the asymptotic behavior of $\alpha _n$'s:
\begin{proposition}[\cite{math.CA/0502300}] \label{cor:verblunskyFiniteCase}
Under assumptions of Theorem \ref{prop:many dominant}, the
Verblunsky coefficients satisfy
\begin{equation*}\label{Verblunsky_asymptotics_dominant_poles}
 \alpha _n=
   -\sum_{k=1 }^\ell \binom{n+1}{m-1}\, \overline{a_k^{n-m+1 }\,
 D_{\rm i} (w; a_k)\, \widehat D_{\rm e} (w; a_k ) }  + \begin{cases}
    \mathcal O \left( \rho ^n \delta ^n\right), & \text{if } \, m=1\,, \\
    \mathcal O\left(  n^{m-2} \, \rho ^n\right)\,, & \text{if }\, m\geq 2\,.
    \end{cases}
\end{equation*}
\end{proposition}

That is, in the situation when the first singularities of $D_{\rm
e}$ met during its analytic continuation inside are only poles, the
Verblunsky coefficients are asymptotically equal  to a combination
of competing exponential functions with coefficients that are
polynomials in $n$. We can compare it with the case of the essential
singularity considered before: for the weight $w$ given in
\eqref{weightExampleSing},
$$
\alpha_n=-\frac{1}{2\sqrt{\pi}}\, t_+^{n} \FF(w; t_+)\, \left(
\frac{\rho }{n}\right)^{3/4} \left( 1+\mathcal O \left(
\frac{1}{n^{1/2}}\right)\right)\,, \quad n \to \infty\,,
$$
where $t_+ \to \rho$ is given by the equation \eqref{equationForT}.

\medskip

A similar analysis can be carried out for the asymptotic expansion
of the leading coefficients $\kappa _n$. By \eqref{kappas},
$$
\kappa_{n}^2=\frac{\tau^2}{2 \pi}\,  S_{22}(n+1; 0)=\frac{\tau^2}{2
\pi}\,\left( 1+g_{n+1}^{(2)}(0)+\mathcal O(r^{4n})\right)\,.
$$
Taking into account \eqref{g1}, we arrive at
\begin{proposition}
For the leading coefficient $\kappa_{n}$ the following formula
holds:
\begin{equation}\label{kappa_new_formula}
\kappa_{n}^2= \frac{\tau^2}{2 \pi}\, \sum_{k> -n-1 } |\left( \FF
\right)_k|^2 +\mathcal O(r^{4n}) =\frac{\tau^2}{2 \pi}\, \left(
1-\sum_{k\leq  -n-1 } |\left( \FF \right)_k|^2 \right)+\mathcal
O(r^{4n})  \,,
\end{equation}
where $\left( \FF \right)_{k}$'s are the coefficients of the Laurent
expansion of $\FF$ in \eqref{LaurentExpansion}.
\end{proposition}
Observe that we can write this identity also in terms of the Riesz
projections:
$$
\kappa_{n}^2=\tau^2 \left\| \mathcal P_-\left( \sigma_{n+1}^{-1}
 \right)\right\|^2_{L^2(\T_1, |dz|)}  +\mathcal O(r^{4n}) \,.
$$
Formula \eqref{kappa_new_formula} shows that
$$
\lim_{n} \kappa_{n}^2=\frac{\tau^2}{2 \pi}\,, \quad \text{and} \quad
\kappa_{n+1}^2-\kappa_{n}^2=\frac{\tau^2}{2 \pi}\, |(\FF)_{-n-1}|^2
+\mathcal O(r^{4n}) \,,
$$
in accordance with \eqref{asympt_expr_alpha} and the well known fact
that
$$
\frac{1}{\kappa _{n+1}^2}-\frac{1}{\kappa _n^2}=-\frac{|\alpha
_n|^2}{\kappa _n^2}.
$$

Summarizing, we see that the Laurent coefficients of $\FF$ (or of
$1/\FF$) contain surprisingly good approximations of two main
parameters of the OPUC: they match asymptotically the Verblunsky
coefficients, and the partial sums of the squares of their absolute
values represent (up to a normalizing constant) the leading
coefficient of the orthonormal polynomials.

\section{Weight with zeros on $\T_1$}\label{sec:nonanalytic}

Let us analyze the change of the behavior of the orthogonal
polynomials if we allow zeros of the weight on the unit circle. In
other words, we consider now a weight of the form
\begin{equation}\label{THE_weight}
W(z) \isdef  w(z)\, \prod_{k=1}^m |z-a_k|^{2\beta _k} \,, \quad z\in
\T_1 \,,
\end{equation}
where $a_k \in \T_1 $, $\beta _k\geq 0$, $k=1, \dots, m$, and $w$ is
an analytic and positive weight on $\T_1$, such as considered in
Section \ref{sec:analytic}. Without loss of generality we assume
that $w$ is analytic and non-vanishing in the annulus $\rho
<|z|<1/\rho $.

According to Nevai and Totik \cite{Nevai/Totik:89}, the Verblunsky
coefficients no longer have an exponential decay, neither the bulk
of zeros accumulate on an inner circle, but how many of them stay
inside? And for those approaching the unit circle, does the rate
depend on the ``orders'' $\beta _k$? And how can we extend the
Riemann-Hilbert method, that so nicely worked for us in the analytic
situation, to the case of a weight of the form \eqref{THE_weight}?

\subsection{Steepest descent analysis}
\label{sec:RH_analysis2}

\begin{figure}[htb]
\centering \begin{overpic}[scale=0.7]{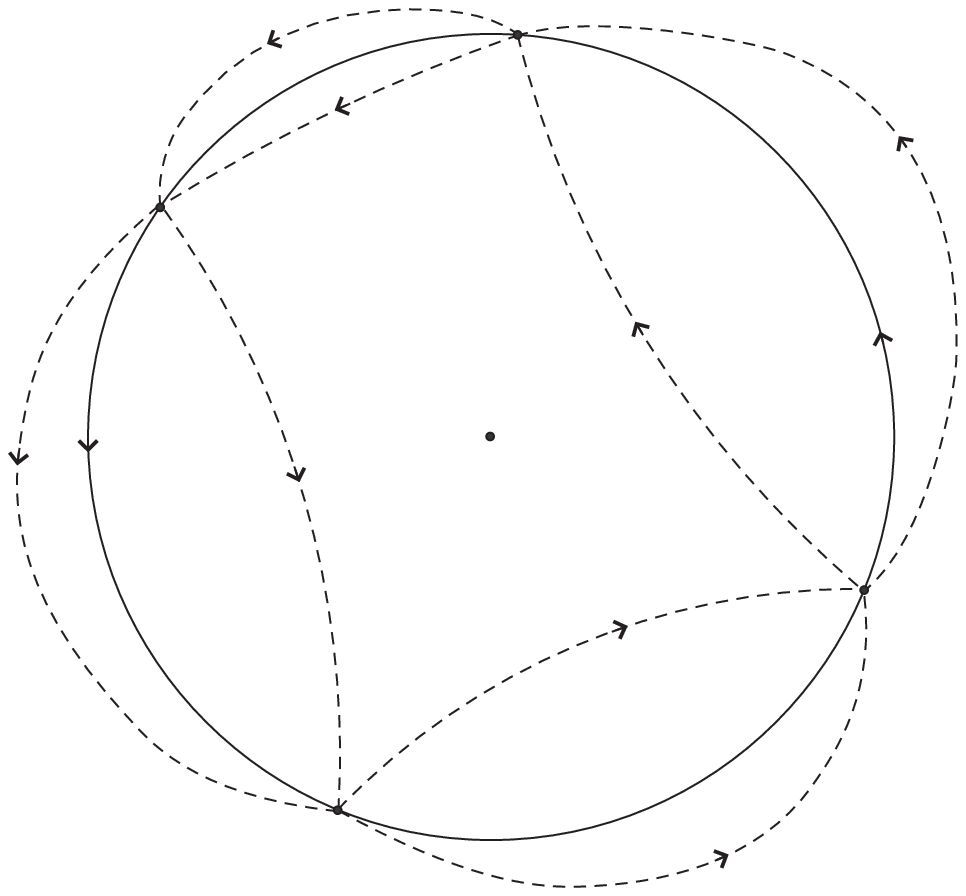}
      \put(60,14){\small $\T_1$}
         \put(36,32){$\gamma_{\rm i}$}
         \put(1,32){$\gamma_{\rm e}$}
          \put(45,65){$\Omega_0$}
          \put(69,69){$\Omega_+$}
          \put(81,76){$\Omega_-$}
          \put(100,70){$\Omega_\infty$}
          \put(22,68){\small $a_m$}
          \put(32,8){\small $a_1$}
         \put(90,31){\small $a_2$}
\end{overpic}
\caption{Opening lenses.}\label{fig:zeros_case1}
\end{figure}

Matrix
$$
Y(z)=\begin{pmatrix} \Phi_n(z) & \displaystyle \dfrac{1}{2\pi
i}\,\oint_{\strut\T_1}
\dfrac{\Phi_n(t) W(t)\, dt}{t^n(t-z)}  \\
-2\pi \kappa _{n-1} \varphi^*_{n-1}(z) & - \, \displaystyle
\dfrac{\kappa _{n-1}}{ i}\,\oint_{\strut\T_1}
\dfrac{\varphi^*_{n-1}(t) W(t)\, dt}{t^n(t-z)}
\end{pmatrix}\,,
$$
solves the Riemann-Hilbert problem \eqref{RHproblem}, with $w$
replaced by $W$. It is the unique solution if we add additional
requirements at the zeros of the weight:
$$
Y(z)= \mathcal O \begin{pmatrix}  1 & 1 \\ 1 & 1 \end{pmatrix} \,,
\quad \text{as } z\to a_k, \; z\in \C\setminus \T_1, \quad k=1,
\dots, m.
$$
Nothing hinders performing Step 1: with $H$ defined by \eqref{defH}
we put $T\isdef Y H$, so that $T$ becomes holomorphic in
$\C\setminus {\T_1}$ (including the infinity) and
$$
T_+(t)=T_-(t)\, \begin{pmatrix} t^n & W(t)  \\ 0 & t^{-n}
\end{pmatrix}\,, \quad t\in \T_1\,.
$$
However, in order to get rid of the oscillatory behavior of the
diagonal entries of the jump matrix the lenses we opened in Step 2
of Section \ref{sec:analytic} are no longer valid, at least because
$W$ has singularities on $\T_1$. Since they are in a finite number,
we can modify this step by opening lenses inside and outside $\T_1$,
but ``attached'' to the unit circle at $a_k$'s (see Figure
\ref{fig:zeros_case1}).
This deformation of the contour makes the definition of $K$ in
\eqref{KforT} consistent (after replacing $w$ with $W$), in such a
way that $U \isdef T K$ has the jumps
$$ U_+(t)=U_-(t)\, J_U(t), \quad
t\in {\gamma_{\rm i}} \cup {\T_1}\cup {\gamma_{\rm e}},
$$
with
\begin{equation}\label{def_JU_for_zeros}
J_U(t)=\begin{cases} \begin{pmatrix} 0 & W(t)  \\ -1/W(t) & 0
\end{pmatrix}, & \text{if } t\in {\T_1}, \\   \begin{pmatrix} 1 & 0  \\ t^n/  W(t)  & 1
\end{pmatrix}, & \text{if } t\in {\gamma_{\rm i}}, \\
\begin{pmatrix} 1 & 0  \\ 1/(t^n W(t)) & 1
\end{pmatrix}, & \text{if } t\in {\gamma_{\rm e}}.
\end{cases}
\end{equation}

{\sc Step 3:} Our next goal is to handle the jump on $\T_1$ via the
global parametrix $N$ built in Section \ref{sec:analytic} using the
Szeg\H{o} function. However, for $W$ as in \eqref{THE_weight} the
Szeg\H{o} function is, in general, no longer single-valued in the
neighborhood of $a_k$'s, and a short digression is convenient in
order to discuss briefly the form of this function and its
multivaluedness.

For the sake of brevity  we define the set of singularities of the
weight, $\mathcal A \isdef \{a_1, \dots, a_m \}$. We fix for what
follows $0 <\delta <1-\rho $, such that additionally $\delta
<\frac{1}{3}\min_{i\neq j} |a_i-a_j|$, so that all neighborhoods
$B_\delta(a_k)$ (see definition \eqref{def_neighborhoods}) are
disjoint. Denote also
\begin{equation}\label{notationBC}
    c_k\isdef \{z\in \C:\, |z-a_k|=\delta\}
    \,, \quad k=1, \dots, m\,,
\end{equation}
as well as $B\isdef \cup_{k=1}^m B_\delta(a_k)$. Furthermore, given
a subset $X\subset \C$ and a value $a\in \C$ we will use the
standard notation $ a \cdot X =\{a x:\, x\in X \}$; consistently,
$\AA \cdot X\isdef \cup_{k=1}^m (a_k\cdot X)$.

In order to construct explicitly the Szeg\H{o} function for the
modified weight $W$ we introduce the generalized polynomial
\begin{equation}\label{def_q}
q(z)\isdef \prod_{k=1}^m (z-a_k)^{\beta _k/2}
\end{equation}
and select its single-valued analytic branch in $\C\setminus
\left(\cup_{k=1}^m a_k\cdot [1,+\infty) \right)$ by fixing the value
of $q(0)$. With this convention we can write the Szeg\H{o} functions
for the modified weight $W$:
\begin{equation}\label{def:D_eandD_i}
D_{\rm i}(W; z)=   \frac{q^2(z)}{q^2(0)} \, D_{\rm i}(w; z)\,, \quad
D_{\rm e}(W; z)=   \frac{D_{\rm e}(w; z)}{q^2(0)\, \overline{q
(1/\bar z)}^2 }\,.
\end{equation}
In particular, $D_{\rm i}(W; z)$ is holomorphic in $\D_{1/\rho}
\setminus \left(\cup_{k=1}^m a_k\cdot [1,1/\rho ) \right)$, $D_{\rm
e}(W; z)$ is holomorphic in $\{ z\in\C:\, |z|>\rho \} \setminus
\left(\cup_{k=1}^m a_k\cdot (\rho,1] \right)$, and
\begin{equation}\label{tauForW}
 \frac{1}{D_{\rm i}(W; 0)}=\frac{1}{ D_{\rm i}(w; 0)}=D_{\rm e}(W;
\infty)=D_{\rm e}(w; \infty)=\tau>0\,,
\end{equation}
where $\tau $ has been defined in \eqref{tau}.  Furthermore, with
the orientation of the cuts toward infinity we have for $ k=1,
\dots, m$:
\begin{equation}\label{boundaryValueD_onCuts}
\begin{split}
\left[D_{\rm i}(W; z) \right]_+ & =e^{-2\pi i \beta _k}\left[D_{\rm
i}(W; z) \right]_-\,, \quad z \in a_k\cdot(1,
    1/\rho )\,,  \\ \left[D_{\rm e}(W; z) \right]_+ & =e^{-2\pi i \beta _k}\left[D_{\rm
e}(W; z) \right]_-\,, \quad z \in a_k\cdot (\rho ,
    1 )\,.
\end{split}
\end{equation}
By definition \eqref{def_F} and formulas \eqref{def:D_eandD_i} we
have
\begin{equation*}
    \FF(W;z) = D_{\rm i}(W; z)D_{\rm e}(W; z)=
  \left( \frac{  q (z) }{q^2(0)\, \overline{q (1/\bar
z)}  } \right)^2 \FF(w; z)   \,,
\end{equation*}
that is also analytic and single-valued in the cut annulus
$\{\rho<|z|<1/\rho \}  \setminus \left(\AA \cdot (\rho ,1/\rho )
\right)$; furthermore, with our assumptions on $\delta $, function
\begin{equation}\label{defSModified}
    \widehat \FF_k(W;z)\isdef \begin{cases}
e^{ \pi i \beta _k} \FF(W; z)\,, & \text{if } z\in B_\delta(a_k)
\text{ and } \arg(z)>\arg(a_k)\,, \\ e^{- \pi i \beta _k} \FF(W;
z)\,, & \text{if } z\in B_\delta(a_k) \text{ and } \arg(z)<
\arg(a_k)\,,
    \end{cases}
\end{equation}
is holomorphic in $B_\delta(a_k)$,  $k=1, \dots, m$. So, we can
define
\begin{equation}\label{def_new_constant}
    \vartheta_k \isdef \widehat \FF_k(W;a_k)\in \T \,, \quad k=1, \dots,
    m\,.
\end{equation}

Now let us get back to the global parametrix $N(z)=N(W; z)$, given
by formula \eqref{equ:defN}, that is well defined, has the same
jumps on $\T_1$ as $U(z)$, and by \eqref{tauForW}, it exhibits the
same behavior at infinity. Hence, $U(z)N^{-1}(z)$ tends to $I$ as
$z\to \infty$, and is holomorphic in $\C\setminus ( \gamma_{\rm
e}\cup \gamma_{\rm i})$. The jump on these curves is again
exponentially close to identity, \emph{except in a neighborhood of
the zeros $a_k$ of the weight $W$.} This is a new feature, and we
have to deal with this problem separately.

\medskip

{\sc Step 4:} local analysis. Let us pick a singular point $a_k\in
\AA$. For the sake of brevity along this subsection we use the
following shortcuts for the notation: $a\isdef a_k$, $\beta \isdef
\beta _k$, $b\isdef B_\delta (a_k)$, $c\isdef c_k$ (where  $c_k$
were defined in \eqref{notationBC}), and $b^+\isdef \{z\in b:\,
\arg(z)>\arg(a) \}$, $b^-\isdef \{z\in b:\, \arg(z)<\arg(a) \}$.  We
also write $\widehat \Omega_j \isdef \Omega_j \cap b$, where $j \in
\{{\rm i}, {\rm e} , 0, \infty \}$, and analogous notation for
curves: $\widehat \T_1 \isdef \T_1  \cap b$, etc.

The goal is to build a matrix $P$ such that it is holomorphic in $b
\setminus (\T_1  \cup \gamma_{\rm i} \cup \gamma_{\rm e})$,
satisfies across $\widehat \T_1  \cup \widehat \gamma_{\rm i} \cup
\widehat \gamma_{\rm e}$ the jump relation $P_+=P_- J_U$, with $J_U$
given in \eqref{def_JU_for_zeros}, with the same local behavior as
$U$ close to $z=a$, and matching $N$ on $c$. This analysis is very
technical, and we refer the reader to \cite{AMF/McL/Saff:06} for
details, and describe here the main ideas in a very informal
fashion.

\begin{figure}[htb]
\centering \begin{overpic}[scale=0.9]{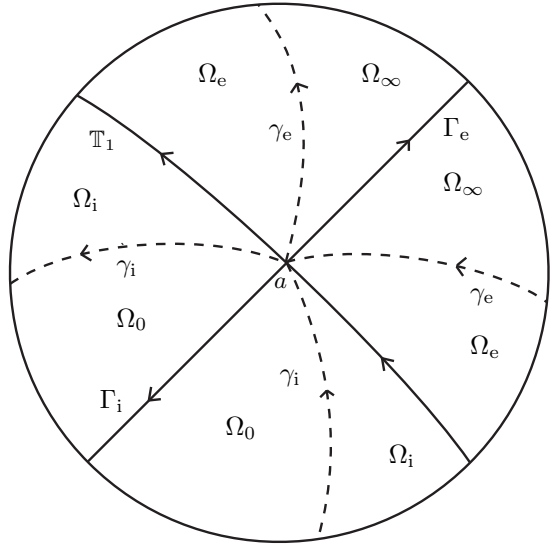}
         \put(50,30){$\gamma_{\rm i}$}
         \put(20,50){$\gamma_{\rm i}$}
         \put(48,75){$\gamma_{\rm e}$}
         \put(85,45){$\gamma_{\rm e}$}
          \put(40,20){$\Omega_0$}
          \put(20,40){$\Omega_0$}
           \put(65,85){$\Omega_\infty$}
          \put(80,65){$\Omega_\infty$}
          \put(35,85){$\Omega_{\rm e}$}
          \put(85,35){$\Omega_{\rm e}$}
          \put(12,62){$\Omega_{\rm i}$}
          \put(70,15){$\Omega_{\rm i}$}
                    \put(17,25){$\Gamma_{\rm i}$}
                   \put(80,75){$\Gamma_{\rm e}$}
         \put(49,47){\small $a$}
       \put(15,73){\small $\T _1$}
\end{overpic}
\caption{Local analysis in $b$.}\label{fig:local1}
\end{figure}

As a first step  we  reduce the problem to the one with constant
jumps. Let us denote $\Gamma_{\rm i}\isdef a\cdot (1-\delta, 1)$ and
$\Gamma_{\rm e}\isdef a\cdot (1, 1+\delta)$, oriented both from
$z=a$ to infinity (see Fig. \ref{fig:local1}). Let $w^{1/2}(z) $ and
$z^{1/2}$ denote the principal holomorphic branches of these
functions in $b$, and $ W^{1/2}(z)\isdef q(z)
\overline{q(1/\overline{z})}\, w^{1/2}(z)$, with $q$ defined in
\eqref{def_q}. Then $W^{1/2}$ is holomorphic in $b\setminus a\cdot
(1-\delta,1+\delta)$, and according to
\eqref{boundaryValueD_onCuts},
\begin{equation*}
\frac{W^{1/2}_+(z)}{W^{1/2}_-(z)}= e^{- \pi i \beta }\quad \text{on
} \Gamma_{\rm i}\,, \quad \text{and} \quad
\frac{W^{1/2}_+(z)}{W^{1/2}_-(z)}= e^{  \pi i \beta }\quad \text{on
} \Gamma_{\rm e}\,.
\end{equation*}
Thus, if we define
\begin{equation}\label{lambda}
\lambda( \beta ;z) \isdef \begin{cases} e^{\pi i \beta } W^{1/2}(z)
z^{n/2}, & z
\in  (\widehat \Omega_{\rm e}  \cup \widehat \Omega_\infty   ) \cap b^+  ,\\
e^{-\pi i \beta } W^{1/2}(z) z^{n/2}, & z \in  (\widehat \Omega_{\rm
e}  \cup \widehat \Omega_\infty   ) \cap b^-  ,\\  e^{-\pi i \beta }
W^{1/2}(z) z^{-n/2}, & z
\in  (\widehat  \Omega_{\rm i}  \cup \widehat \Omega_0   ) \cap b^+  ,\\
e^{\pi i \beta } W^{1/2}(z) z^{-n/2}, & z \in  (\widehat \Omega_{\rm
i}  \cup \widehat \Omega_0   ) \cap b^-  ,
\end{cases}\,,
\end{equation}
and set
\begin{equation}\label{def_R}
R(z) \isdef  P(z)\, \lambda (  \beta ; z)^{  \sigma _3} \,, \quad z
\in b \setminus (\Gamma_{\rm i} \cup \Gamma_{\rm e} \cup \T _1 \cup
\gamma_{\rm i} \cup \gamma_{\rm e})\,,
\end{equation}
we get for $R$ the following problem: $R$ is holomorphic in $b
\setminus (\Gamma_{\rm i} \cup \Gamma_{\rm e} \cup \T_!  \cup
\gamma_{\rm i} \cup \gamma_{\rm e})$, and satisfies the jump
relation $R_+(z)=R_-(z) J_R(z)$, with
\begin{equation*}
J_R(z)=\begin{cases}  \begin{pmatrix} 0 & 1  \\ -1  & 0
\end{pmatrix}, & \text{if } z \in \widehat \T_1   , \\
 \begin{pmatrix} 1 & 0  \\  e^{-2\pi i \beta  }  & 1
\end{pmatrix}, & \text{if } z\in   (\widehat \gamma_{\rm i}   \cap b^+)\cup
(\widehat \gamma_{\rm e}   \cap b^-)  ,
\\
\begin{pmatrix} 1 & 0  \\  e^{ 2\pi i \beta  }  & 1
\end{pmatrix}, & \text{if } z\in   (\widehat \gamma_{\rm i}   \cap b^-)\cup (\widehat \gamma_{\rm e}   \cap b^+) ,\\
\begin{pmatrix} e^{ \pi i \beta } & 0  \\  0  & e^{-\pi i \beta }
\end{pmatrix}, & \text{if } z\in \Gamma_{\rm i}\cup \Gamma_{\rm e}\,.
\end{cases}
\end{equation*}
Moreover, $R$ has the following local behavior as $z\to a$:
\begin{equation*}
R(z)=\begin{cases} \mathcal O \begin{pmatrix}  |z-a|^{\beta  } & |z-a|^{-\beta  } \\
|z-a|^{\beta  } &
|z-a|^{-\beta  } \end{pmatrix}, & \text{if } z \in \widehat \Omega_0 \cup \widehat \Omega_\infty, \\
\mathcal O \begin{pmatrix}  |z-a|^{-\beta  } & |z-a|^{-\beta  } \\
|z-a|^{-\beta  } & |z-a|^{-\beta  }
\end{pmatrix}, & \text{if } z \in \widehat \Omega_{\rm e} \cup \widehat \Omega_{\rm
i}.
\end{cases}
\end{equation*}

Consider in $\C \setminus (-\infty, 0)$ the transformation
\begin{equation}\label{def_zeta_local_analysis}
\zeta = - i\, \frac{n}{2}\, \log(z/a)\,,
\end{equation}
(we omit the explicit reference to the dependence of $\zeta $ from
$a$ and $n$ in the notation), where we take the main branch of the
logarithm. This is a conformal 1-1 map of $b$ onto a neighborhood of
the origin. Moreover, $\T_1 $ is mapped onto $\R$ oriented
positively, $\Gamma_{\rm i} \cup \Gamma_{\rm e}$ are mapped on the
imaginary axis, and we may use the freedom in the selection of the
contours deforming them in such a way that $f(\widehat \gamma_{\rm
i})$ and $f(\widehat \gamma_{\rm e})$ follow the rays $\{\arg \zeta
=\pm \frac{\pi}{4}\pm \pi \}$. After this transformation we get a
Riemann-Hilbert problem on the $\zeta$-plane that has been studied
for the local analysis of the generalized Jacobi weight on the real
line. We take advantage of the results proved therein in order to
abbreviate the exposition, and refer the reader to  \cite[Theorem
4.2]{Vanlessen03} where the solution $ \Psi \left( \beta ; \zeta
\right)$ is explicitly written in terms of the Hankel  and modified
Bessel functions.

Since a left multiplication by a holomorphic function has no
influence on the jumps, and taking into account \eqref{def_R}, we
see that matrix $P$ can be built of the form
\begin{equation}\label{formOfP}
P(z) = E(z)\, \Psi \left( \beta ;  \zeta   \right) \, \lambda (
\beta ; z)^{ - \sigma _3}\,,
\end{equation}
where $E$ is any holomorphic function in $b$. An adequate selection
of $E$ is motivated by the matching requirement $P(z) N^{-1}(z
)=I+\mathcal O (n^{-1})$ on the boundary $c$, and is constructed
analyzing the asymptotic behavior of the matrix-valued function
$\Psi \left( \beta ;  \zeta   \right)$ at infinity. Let us summarize
the results of this analysis:
\begin{proposition}
\label{local_lemma1} Let
\begin{equation}\label{representation for E}
    E(z)\isdef \left( \frac{\widehat S_k(W;z)}{\tau ^2}\, i a^n \right)^{\sigma _3/2}\frac{1}{\sqrt 2}\,
    \begin{pmatrix}
        i & 1 \\
        -1 & -i
    \end{pmatrix}\,,
\end{equation}
where $\widehat S_k(W;\cdot)$ has been defined in
\eqref{defSModified}, and we take the main branch of the square
root. Then matrix $P(z)=P(a,\beta ;z)$,
\begin{equation}\label{def_P_for_localA}
    P(a,\beta ;z)\isdef E(z)\, \Psi \left( \beta ;  \zeta  \right) \,
\lambda (  \beta ; z)^{ - \sigma _3}\,,
\end{equation}
with $\zeta $ given by \eqref{def_zeta_local_analysis}, satisfies:
\begin{enumerate}
\item[(i)] $U(z)P^{-1}(z)$ is holomorphic in $b$;
\item[(ii)] for $z \in c$,
\begin{equation}\label{matching_cond}
\begin{split}
 P(z) N^{-1}(z) = & I+\frac{i \beta   }{2 \zeta }   \,
\begin{pmatrix}
             \beta      &      - \tau ^{-2}    a^n \widehat \FF
    (W;z)       \\[2ex]
          \tau ^{2}      a^{-n}\widehat \FF^{-1}
    (W;z)   &  - \beta
    \end{pmatrix}   +  \mathcal O\left(\frac{1}{n^2}
    \right)\,.
\end{split}
\end{equation}
In particular, $P(z) N^{-1}(z)=I+\mathcal O (n^{-1})$ for $z\in c$.
\end{enumerate}
\end{proposition}

\medskip

{\sc Step 5:} asymptotic analysis. With the notation introduced in
\eqref{notationBC} and with $P(a,\beta ;z)$ defined by
\eqref{representation for E}--\eqref{def_P_for_localA} let us take
\begin{equation*}
P(z)\isdef P (a_k, \beta _k; z )\quad  \text{for } z \in B_\delta
(a_k) \setminus (\T_1  \cup \gamma_{\rm e}\cup \gamma_{\rm i}),
\quad k =1, \dots, m\,,
\end{equation*}
and put
$$
S(z) \isdef \begin{cases} U(z)N^{-1}(z), & \text{for } z \in
\C\setminus (B \cup  \T_1  \cup \gamma_{\rm e}\cup \gamma_{\rm i}), \\
U(z)P^{-1}(z), & \text{for } z \in B  \setminus (\T_1  \cup
\gamma_{\rm e}\cup \gamma_{\rm i}).
\end{cases}
$$
\begin{figure}[htb]
\centering \begin{overpic}[scale=0.7]{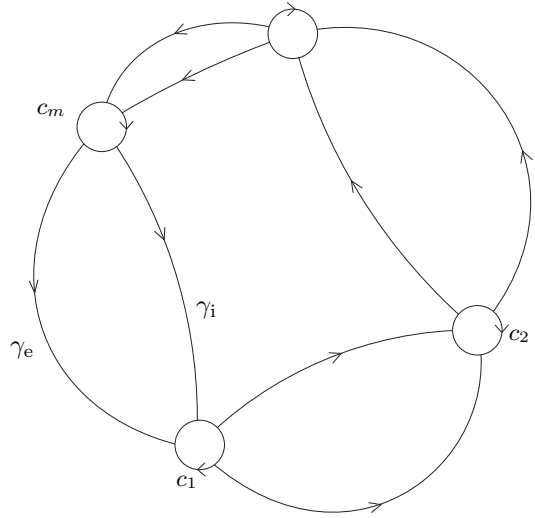}
         \put(35,45){$\gamma_{\rm i}$}
         \put(7,39){$\gamma_{\rm e}$}
          \put(11,75){\small $c_m$}
          \put(32,19){\small $c_1$}
         \put(82,41){\small $c_2$}
\end{overpic}
\caption{Jumps of $S$.}\label{fig:lensesFinal}
\end{figure}
Matrix $S$ is holomorphic in the whole plane cut along $  \gamma
\cup C$, where
$$
\gamma  \isdef (\gamma_{\rm e}\cup \gamma_{\rm i}) \setminus B \quad
\text{and} \quad C\isdef \cup_{k=1}^m c_k
$$
(see Fig. \ref{fig:lensesFinal}), $S(z) \to I$ as $z\to \infty$, and
if we orient all $c_k$'s clockwise, $S_+(t)=S_-(t) J_S$, with
$$
J_S(t)=\begin{cases} P (z) N^{-1}(z), & \text{if } z\in C,  \\
\begin{pmatrix}
1 & 0 \\ \tau^2 /(z^n \FF(W;z)) & 1
\end{pmatrix}, & \text{if } z\in \gamma_{\rm e}\setminus B, \\
\begin{pmatrix}
1 &  - z^n \FF(W; z)/\tau^2  \\ 0 & 1
\end{pmatrix}, & \text{if } z\in \gamma_{\rm i}\setminus B. \\
\end{cases}
$$
It is clear that the off-diagonal terms of $J_S$ on $\gamma_{\rm
i}\setminus B$ and $\gamma_{\rm e}\setminus B$ decay exponentially
fast. On the other hand, by \eqref{matching_cond},
$J_S(z)=I+\mathcal O(1/n)$ for $z\in C$. So the conclusion is that
the jump matrix $J_S=I + O(1/n)$ uniformly for $z \in  \gamma \cup C
$. Then arguments such as in
\cite{MR2000g:47048,MR2001f:42037,MR2001g:42050} lead to the
following conclusion:
\begin{proposition}
Matrix $S$ satisfies the following singular integral equation:
\begin{equation*}
S(z)=I + \frac{1}{2\pi i}\, \int \frac{ (J_S(t)-I)\, dt}{t-z} +
\frac{1}{2\pi i}\, \int  \frac{ (S_- (t)-I)(J_S(t)-I)\, dt}{t-z}\,,
\end{equation*}
where we integrate along contours $  \gamma \cup C $ with the
orientation shown in Fig.\ \ref{fig:lensesFinal}. In particular,
integrating counter-clockwise,
\begin{equation*} 
    S(z) = I - \sum_{k=1}^m \frac{1}{2\pi i}\, \oint_{c_k} \frac{ (P(t) N^{-1}(t)-I)\,
dt}{t-z} +   O\left(\frac{1}{n^2}\right)
    \end{equation*}
locally uniformly for $z \in \C\setminus ( \gamma \cup C).$
\end{proposition}

Now formula \eqref{matching_cond} and the residue theorem yield for
$z \in \C\setminus (\gamma\cup B)$,
\begin{equation}\label{residues_Final}
S(z) = I + \frac{1}{n}\, \sum_{k=1}^m  \dfrac{a_k \beta _k  }{
 a_k-z } \,
\begin{pmatrix}
             \beta_k      &      - \tau ^{-2} a_k^n \vartheta_k    \\[2ex]
         \tau ^{ 2}  a_k^{-n} \vartheta_k^{-1}  &  - \beta_k
    \end{pmatrix} +   O\left(\frac{1}{n^2}\right)\,.
\end{equation}
We are ready for the asymptotic analysis of the original matrix $Y$
(and in particular, of its entries $(1,1)$ and $(2,1)$).

\begin{figure}[htb]
\centering \begin{overpic}[scale=0.7]{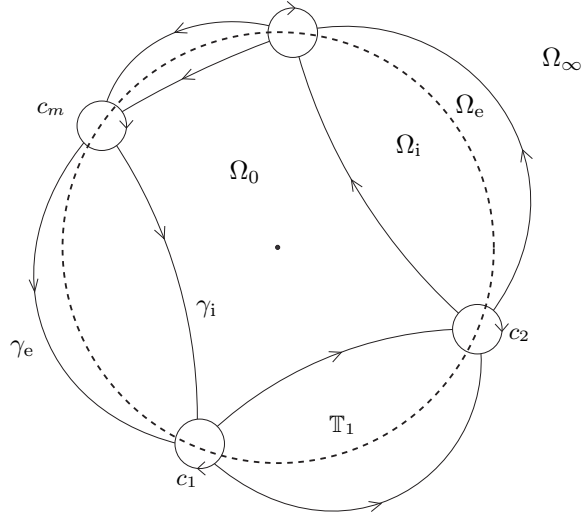}
      \put(55,27){\small $\T_1 $}
         \put(35,45){$\gamma_{\rm i}$}
         \put(7,39){$\gamma_{\rm e}$}
          \put(40,65){$\Omega_0$}
          \put(65,69){$\Omega_{\rm i}$}
          \put(74,75){$\Omega_{\rm e}$}
          \put(87,82){$\Omega_\infty$}
          \put(11,75){\small $c_m$}
          \put(32,19){\small $c_1$}
         \put(82,41){\small $c_2$}
\end{overpic}
\caption{Domains for the asymptotic
analysis.}\label{fig:lensesFinalAsymptotics}
\end{figure}

Unraveling our transformations we have
\begin{equation*}
    Y(z)=\begin{cases}
     S(z)N(z)K^{-1}(z)\, H^{-1}(z), & \text{if } z \in \C \setminus B , \\  S(z)
P(z)K^{-1}(z)\, H^{-1}(z), & \text{if } z \in  B.
    \end{cases}
\end{equation*}
We must analyze the consequences of these formulas in each domain
(see Fig.\ \ref{fig:lensesFinalAsymptotics}).

We will do it only for the interior domain $\Omega_0\setminus B$,
where
$$
N(z)=    \begin{pmatrix} 0 & D_{\rm i}(W; z)/\tau   \\ -\tau /D_{\rm
i}(W; z) & 0
\end{pmatrix}, \quad K(z)=H(z)=I\,.
$$
Hence, $  Y(z)=S(z) N(z)$, so that
$$
Y_{11}(z)=-\frac{\tau }{D_{\rm i}(W; z)}\, S_{12}(z)\,, \quad
Y_{21}(z)= -\frac{\tau }{D_{\rm i}(W; z)}\, S_{22}(z)\,.
$$
Taking into account \eqref{tauForW} and \eqref{residues_Final}, and
recalling that $\Phi_n =Y_{11} $ we obtain
\begin{equation*}
\Phi_n(z)=  \frac{ D_{\rm i}(W; 0)}{  D_{\rm i}(W; z)}\,
\frac{1}{n}\, \left( \sum_{k=1}^m \dfrac{\beta _k \vartheta_k
 }{ a_k-z }\, a_k^{n+1}    + O\left(\frac{1}{n }\right) \right)\,,
\end{equation*}
valid uniformly in this domains. It shows in particular that for
every compact set $K\subset \D_1$ there exists $N=N(K)\in \N$ such
that for every $n\geq N$, each $\Phi_n$ has at most $m-1$ zeros on
$K$, and these zeros should be asymptotically close to those of the
rational fraction
\begin{equation} \label{where_zeros}
 \sum_{k=1}^m \dfrac{\beta _k \vartheta_k
 }{ a_k-z }\, a_k^{n+1} \,.
\end{equation}
Evaluating $\Phi_{n}(z)$ at $z=0$ we can obtain asymptotics for the
Verblunsky coefficients $\alpha _n$; I leave this as an exercise for
an interested reader.

\medskip

{\sc Remark:} Orthogonal polynomials with respect to non-analytic
(but smooth) and non-vanishing weights with finite jump
discontinuities on $\T_1$ were studied in
\cite{McLaughlin/Miller:2004}, using a different (but complementary)
method based on the $\overline
\partial$ problem. There are very interesting similarities
with our case: roughly speaking, both zeros and jump discontinuities
of the weight have the same effect on the asymptotics of $\Phi_n$'s
inside the unit disk $\D_1$. More precisely, according to
\cite{McLaughlin/Miller:2004}, if $W>0$ is sufficiently smooth and
piecewise continuous, with jump discontinuities at $a_k\in \T_1$,
$k=1, \dots, m$, then with a proper normalization, $\Phi_n(z)$ for
$z$ on compact subsets of $\D_1$ will be asymptotically close to
rational fractions of the form \eqref{where_zeros}, with the basic
difference that now coefficients $\beta _k$ stand for the magnitude
of the jump of $\log W$ at $a_k$ (cf.\ formula (56) in
\cite{McLaughlin/Miller:2004} and the fact that for the scattering
function on $\T_1$,
$$
\mathcal S\left(W;e^{i\theta}\right)=e^{\Omega(\theta)}\,,
$$
where $\Omega$ is defined by (21) in \cite{McLaughlin/Miller:2004}).
These similarities have as a common ground the duality of both
cases: for $W$ in \eqref{THE_weight}, the imaginary part of $\log W$
has jumps proportional to $\beta_k$'s, while in
\cite{McLaughlin/Miller:2004} the finite jumps correspond to
$\Re\left(\log W\right)$.

\medskip

Finally, recall that  the leading coefficient $\kappa_n$ of the
orthonormal polynomial $\varphi_n$ is expressed in terms of $Y$ by
$Y_{21}(0)=-2\pi \kappa _{n-1}^2$. This yields immediately the
following asymptotic formula:
$$
\kappa _{n-1}^2= \frac{\tau ^2}{2\pi}\,\left( 1 - \frac{1}{n}\,
\sum_{k=1}^m \beta _k^2+\mathcal O\left( \frac{1}{n^2}\right)
\right) \,, \quad n\to \infty.
$$
This result has a consequence for the behavior of the Toeplitz
determinants for $W$. If we define the moments
$$
d_k\isdef \oint_{z\in \T } z^{-n} W(z) |dz|\,,
$$
then the Toeplitz determinants are
\begin{equation}\label{def_Toeplitz_dets}
\mathcal D_n(W)\isdef \det \left[ \left(
d_{j-i}\right)_{i,j=0}^n\right]\,.
\end{equation}
It is known (see e.g.\ \cite[Theorem 1.5.11]{Simon05b}) that
$$
\frac{\mathcal D_n(W)}{\mathcal D_{n-1}(W)}=\frac{1}{\kappa _n^2}\,.
$$
Taking into account the asymptotics of $\kappa _n$, \eqref{def_G}
and \eqref{tauForW}, we arrive at
\begin{theorem} \label{thm:Fisher}
Under the assumption above there exists a constant $\varkappa$
depending on $W$ such that
\begin{equation}\label{asympt_Tplitz}
\begin{split}
\mathcal D_n(W) & =\varkappa\,  \left( G[2\pi w] \right)^n\,
  n^{\sum_{k=1}^m \beta _k^2} \, \left( 1+o (1)\right) \,, \quad n \to \infty\,.
\end{split}
\end{equation}
\end{theorem}
This formula is in accordance with the well known Fisher-Hartwig
conjecture (see e.g.\ \cite{Basor91}), proved for this case (but
using totally different approach and giving an expression for
$\varkappa$) by Widom in \cite{Widom:73}.

\section*{Acknowledgements}

The research of this author was supported, in part, by a grant from
the Ministry of Education and Science of Spain, project code
MTM2005-08648-C02-01, by Junta de Andaluc\'{\i}a, Grupo de Investigaci\'{o}n
FQM229, by ``Research Network on Constructive Complex Approximation
(NeCCA)'', INTAS 03-51-6637, and by NATO Collaborative Linkage Grant
``Orthogonal Polynomials: Theory, Applications and
Generalizations'', ref. PST.CLG.979738.

I wish to acknowledge also the contribution of both anonymous
referees, whose careful reading of the manuscript and kind
suggestions helped to improve the text.

A vast portion of the results exposed here is an outcome of a highly
enjoyable collaboration with K. T.-R. McLaughlin, from the
University of Arizona, USA, and with E.\ B.\ Saff, from the
Vanderbilt University, USA. In particular, this celebration is a
good occasion to express once again my gratitude to Ed Saff for his
support and friendship.

%

\end{document}